\documentclass[11pt,letterpaper]{article}

\usepackage{amsmath,amsthm,amsfonts,amssymb}
\usepackage{mathrsfs}
\usepackage{xcolor}
\usepackage{graphicx}
\usepackage[numbers,sort&compress]{natbib}
\usepackage{subcaption}
\usepackage{algorithm}
\usepackage{algpseudocode}
\usepackage{booktabs}
\usepackage[shortlabels]{enumitem}
\usepackage[letterpaper,margin=1in]{geometry}
\usepackage[hidelinks]{hyperref}
\hypersetup{
  pdftitle={Sub-Riemannian flow matching},
  pdfauthor={Wonjun Lee and Scott Zimmerman}
}

\renewcommand{\epsilon}{\varepsilon}

\def\XXint#1#2#3{{\setbox0=\hbox{$#1{#2#3}{\int}$}
     \vcenter{\hbox{$#2#3$}}\kern-.5\wd0}}

\newtheorem{theorem}{Theorem}[section]
\newtheorem*{theorem*}{Theorem}
\newtheorem{proposition}[theorem]{Proposition}

\theoremstyle{definition}
\newtheorem{definition}[theorem]{Definition}

\theoremstyle{remark}
\newtheorem{remark}[theorem]{Remark}

\numberwithin{equation}{section}

\title{Sub-Riemannian flow matching}
\author{Wonjun Lee\thanks{Corresponding author. \texttt{lee.8222@osu.edu}}
  \quad Scott Zimmerman\thanks{\texttt{zimmerman.416@osu.edu}}\\[4pt]
  \small Department of Mathematics, The Ohio State University, USA}
\date{}

\begin{document}
\maketitle

\begin{abstract}
We study sub-Riemannian conditional flow matching (SRCFM) for generative
transport with prescribed velocity constraints. Training uses paths between
sampled source and target points whose velocities lie in the allowed
horizontal directions at each state. A neural network takes time and the full
state as input and predicts the coefficients of these directions, and the
horizontal frame then determines the full velocity used for generation. We
describe constructions based on Heisenberg geodesics and joint endpoint
equations in free step-two Carnot groups. We propose a scalable piecewise
algorithm that computes horizontal paths directly from their endpoints. A
straight segment sets the horizontal coordinates, and a loop in those
coordinates corrects the remaining vertical displacement. We verify the
algorithm's endpoint accuracy to numerical precision in high-dimensional
tests. Synthetic experiments illustrate approximate distribution transport
along horizontal trajectories. On resized MNIST, we use the same construction
to generate images from Gaussian inputs with fewer controls than pixels.
\end{abstract}

\begin{center}
\end{center}

\section{Introduction}

{
Conditional flow matching (CFM) is a method for generating samples from an unknown probability distribution $\nu$ on $\mathbb{R}^d$.
This is accomplished by
learning a time-dependent velocity field from sampled interpolation
paths and solving its ODE~\cite{lipman2023flowmatching}. 
More precisely, samples are first selected from both $\nu$ and a known, tractable source distribution $\mu$. Second, paths are drawn between endpoints in each sample, and, third, a time dependent vector field $v_\theta$ is learned by averaging the velocities of paths from several samples. This vector field defines an ODE that can be solved to generate a sample from $\nu$ using a sample from $\mu$ by matching endpoints of the flows of $v_\theta$.
Put simply, we attempt to understand a distribution $\nu$ by averaging over samples; CFM does this by averaging paths drawn from these samples.

In standard CFM, a common choice is to join sampled endpoints by straight-line segments. These paths are easy to compute, and flow-matching models based on them have been used successfully for large-scale image generation, unpaired image translation, and modeling single-cell dynamics~\cite{lipman2023flowmatching,liu2023rectifiedflow,tong2024improving}.  In this paper, we consider a $d$-dimensional state space in which the velocity at each point is determined by $m<d$ controls acting along prescribed directions. These directions can vary with the state. We use sub-Riemannian geometry to construct training paths that connect sampled endpoints under this constraint. The constructions include an explicit piecewise path designed to scale to high-dimensional generation.

Fix $m < d$, and suppose $v_1,\dots,v_m$ are vector fields on a state space $\mathbb{R}^d$ that are linearly independent at every point. 
Denote by $F(x)$ the frame of these vector fields, i.e., $F(x)$ is the $d \times m$ matrix whose $j$th column is $v_j(x)$.
We consider paths satisfying the following:
\[
\dot \gamma(t) = F(\gamma(t)) a(t)
\quad
\text{for almost every } t.
\]
In other words, these are paths whose velocities lie in a prescribed $m$ dimensional space at every point. $F$ is the transformation that generates the velocity $\dot \gamma$ from the control data $a$.
We ensure that arbitrary points can be connected by such paths by requiring that the vector fields $v_j$ satisfy the so-called Hormander Condition.

Unlike a constraint requiring that all
states lie on a lower-dimensional surface, our constraint
restricts only the velocity of the paths. 
Once these paths are drawn, the learning of $v_\theta$ in the third step of CFM occurs in the $m$-dimensional control space rather than the full $d$-dimensional state space.

Declaring the vectors $v_j(x)$ to be
orthonormal gives $\mathbb{R}^d$ the structure of a sub-Riemannian manifold. 
The simplest, non-Euclidean sub-Riemannian manifold is the 
the Heisenberg group. 
We will view this as $\mathbb{R}^3$ where admissible curves $\gamma=(x,y,z)$ are those satisfying
\[
 \dot z=-x\dot y.
\]
Only the data for $x$ and $y$ are specified independently;
the velocity of the third component is determined by the above identity.
Note, for example, that vertical segments are not allowed in this space, but any two points in $\mathbb{R}^3$ can be connected by a path satisfying this 2-dimensional control.
More
generally, we will consider direct sums of Heisenberg groups and 
the so-called free Carnot groups of step 2.
The latter case considers an $m$-dimensional control space and a state space of dimension $d=m+\binom m2$.

We study \emph{sub-Riemannian conditional flow matching} (SRCFM).
In this setting, a neural network receives the full state
and time data as in the CFM setting. It then predicts an $m$-dimensional control vector field $a_\theta$ and obtains the vector field $v_\theta$ in the ambient space using the known frame $F$.
This is a reduction of the number of independent velocity outputs rather than the dimension of the state. 
}

The main contribution of this paper is the construction and computational assessment
of SRCFM training paths in sub-Riemannian manifolds of step two. The focus is on the construction of horizontal paths and their use in predicting
regression targets rather than on a new principle of control prediction
or conditional averaging. More specifically:
\begin{enumerate}
 \item We describe minimizing conditional targets using geodesics on Heisenberg groups
 and their direct sums, not-necessarily-geodesic solutions to the geodesic equation on the free step two Carnot groups $\mathbb{R}^d$, and piecewise defined paths on these spaces 
 that can be explicitly defined from the spectral decomposition of a skew-symmetric matrix.
 This last construction avoids the use of a nonlinear endpoint solve required in the previous formulation. 
 \item We place these path choices in a two-endpoint control-regression
 algorithm that reconstructs the ambient sampling velocity from learned horizontal controls.
 \item We use straight-line CFM as an unconstrained reference in the
 Heisenberg experiments, reporting endpoint distributional error and
 horizontal residuals separately. We also test the two path constructions on $\mathbb{R}^d$ with the free step two group structure and illustrate learned
 horizontal transport in an image geometry with selected vertical coordinates.
 These tests do not establish reliable endpoint solving for every value of $m$
 or accurate image-distribution matching.
\end{enumerate}

The significance of our discussion is structural: SRCFM constructs a generation field
whose velocities use the prescribed horizontal directions while
inducing motion in the full ambient space. Euclidean CFM does not
impose this restriction. The comparison therefore asks whether useful
endpoint approximation can be retained while enforcing admissible
motion, not whether SRCFM improves unconstrained generation.
Distributional error and horizontal residual measure different parts
of this requirement. In the Heisenberg tests, Euclidean CFM gives
somewhat smaller mean ambient Wasserstein error, whereas SRCFM gives
much smaller horizontal residuals. The free-group and image studies
also expose remaining distributional errors. Thus horizontality of
the learned field does not by itself establish an accurate
pushforward to the target law. Robust geodesic computation for arbitrary endpoint
pairs in free groups with many horizontal controls remains unresolved in our current
implementation.

\subsection{Related work}
\label{sec:related-work}

\paragraph{Sub-Riemannian geometry and geodesics.}
For a historical reference on the study of free, step 2 Carnot groups using skew-symmetric matrices, see \cite{brockett1980}.
The surveys by Agrachev, Barilari, and Boscain~\cite{agrachev2020} (specifically Chapter 13) and Le Donne \cite{ledonne2017,ledonne2025} are excellent learning resources regarding the geometry and analysis motivating the discussion below.
Foundational results specifically about the particular sub-Riemannian structures discussed below can be found in papers by Myasnichenko \cite{Myasnichenko2006}, Montanari and Morbidelli \cite{montamorbid2017}, and Rizzi and Serres \cite{rizziserres2017}.
More recently, the works of Murray and Sastry~\cite{murray1993steering}, Le Donne, Leonardi, Monti, and Vittone~\cite{ledonne2013extremal}, Haj{\l}asz and
Zimmerman~\cite{hajlaszZimmerman2015}, and Montanari and Morbidelli \cite{montamorbid2026} explore geodesic equations and other curves in the free step 2 Carnot groups.

\paragraph{Flow matching and the choice of paths.}
Flow matching trains continuous-time generative models by regressing velocities along prescribed probability paths~\cite{lipman2023flowmatching}.  Rectified flow uses straight segments between sampled endpoints~\cite{liu2023rectifiedflow}. Conditional and multisample flow matching use different rules for pairing endpoints, including minibatch optimal transport~\cite{tong2024improving,pooladian2023multisample}.  Stochastic interpolants provide a common construction for deterministic flows and diffusions~\cite{albergo2025interpolants}. There are also methods that combine score and flow matching to learn Schr{\"o}dinger bridges~\cite{tong2024schrodinger}.

Geometry enters the choice of training paths in several ways.  Flow matching on Riemannian manifolds uses geodesics or paths computed from distance-like functions~\cite{chen2024flow}. On Lie groups, there is a formulation based on exponential curves~\cite{sherry2026lie}. Related work on neural optimal transport computes least-action paths for prescribed Lagrangian costs~\cite{pooladian2024lagrangian}.

\paragraph{Generation with prescribed dynamics.}
Flow matching has been used to learn controls for stochastic linear systems and nonlinear control-affine systems~\cite{mei2025flow,elamvazhuthi2025flow}. A driftless Heisenberg example constructs training paths by first moving the horizontal coordinates and then correcting the vertical coordinate with sinusoidal controls~\cite[Section~6.3]{elamvazhuthi2025flow}. There are also continuous normalizing flows that learn coefficients of fixed bracket-generating vector fields by maximum likelihood~\cite[Algorithm~1]{morteza2026parameter}.

\paragraph{Sub-Riemannian transport and stochastic bridges.}
There are works on optimal transport with sub-Riemannian distance costs and on dynamical transport for control-affine systems~\cite{figalliRifford2010,elamvazhuthi2023dynamical}. Schr{\"o}dinger bridges with noise along the horizontal directions recover sub-Riemannian optimal transport in the small-noise limit~\cite{adu2026schrodinger}. Score matching has been used to sample endpoint-conditioned sub-Riemannian diffusions, including Heisenberg bridges~\cite{grong2024score}. Mixtures of diffusion bridges have also been used for generative modeling on Riemannian manifolds~\cite{jo2024riemannianmixture}.

We use a specified endpoint coupling and study how to compute deterministic horizontal paths for its sampled pairs. Our focus is the construction and numerical assessment of these paths as training targets in step-two geometries.

\section{Sub-Riemannian geometry}

\subsection{Intuition}
Consider an object moving and rotating in the plane $\mathbb{R}^2$. We will use three variables to describe the state of the object at any time: two coordinates $(u,v)$ to represent its location, and one coordinate $\theta$ to represent its angle with respect to a $u$-axis direction.
That is, the position and orientation of the object in the plane can be described using a point $(\theta,u,v)$ in $\mathbb{R}^3$.
The movement and rotation of the object in the plane can then be viewed as motion through $\mathbb{R}^3$ using the vector fields
\[
V_1(\theta,u,v) = \frac{\partial}{\partial \theta},
\qquad
V_2(\theta,u,v) = \cos(\theta) \frac{\partial}{\partial u} + \sin(\theta) \frac{\partial}{\partial v},
\qquad
V_3(\theta,u,v) = -\sin(\theta) \frac{\partial}{\partial u} + \cos(\theta) \frac{\partial}{\partial v}.
\]
$V_1$ controls the rotation of the object while $V_2$ and $V_3$ control movement in the directions parallel and orthogonal to its orientation respectively.

\begin{figure}
    \centering
    \includegraphics[width=0.5\linewidth]{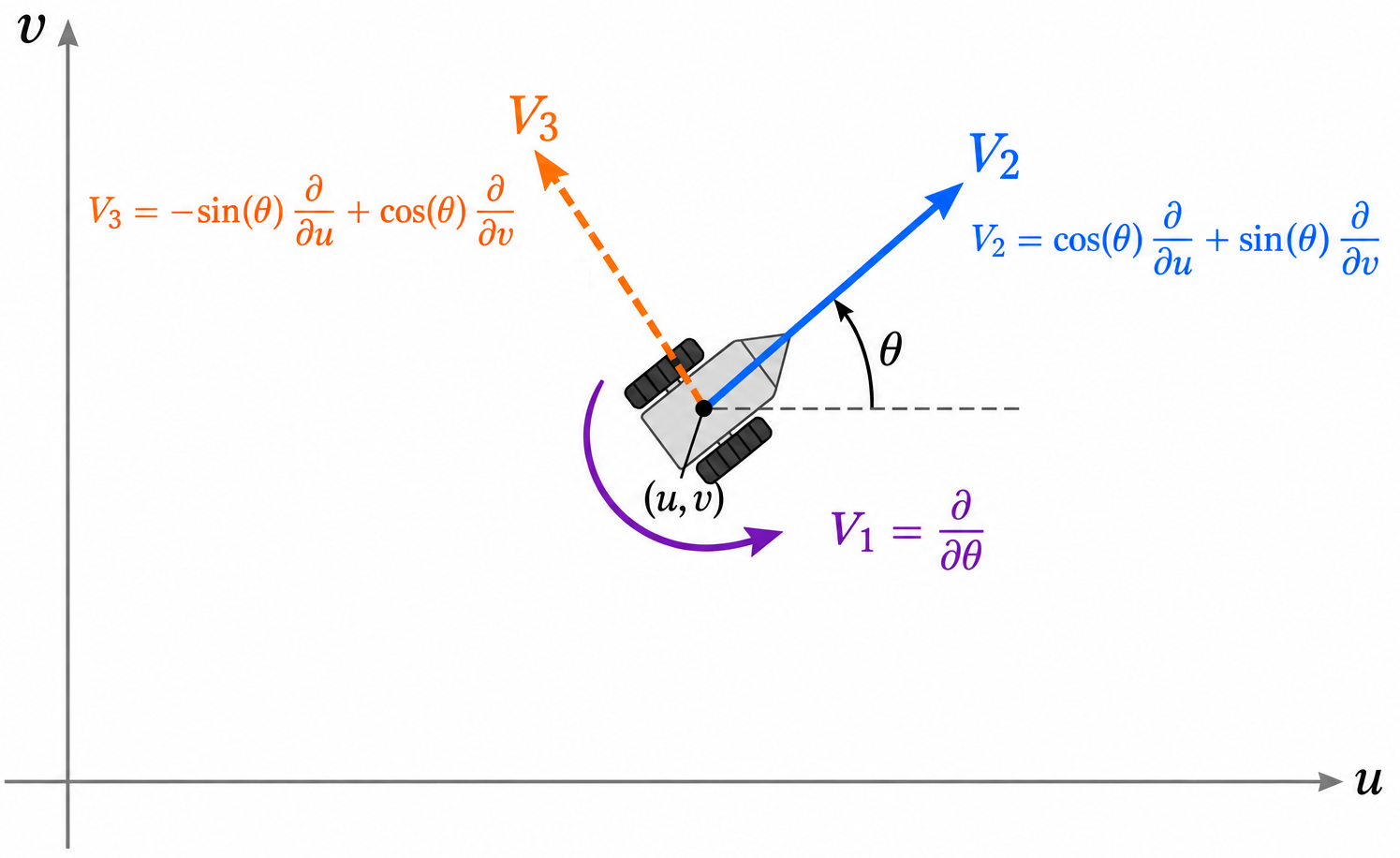}
    \caption{Motion of a two wheeled robot in the $uv$-plane}
    \label{fig:robot}
\end{figure}

Suppose now that the object can no longer move in the $V_3$ direction, such as a robot with two wheels. 
The robot is allowed to drive forward and backward and it can rotate, but it cannot move orthogonal to the direction it is facing, i.e., it cannot move in the $V_3$ direction. See Figure~\ref{fig:robot}. In other words, the motion of the robot in the state space $\mathbb{R}^3$ is restricted at any point $p=(\theta,u,v)$ to the plane of directions
\[
H(p) := \operatorname{span}\{V_1(p),V_2(p)\}.
\]

Observe that, despite the restrictions on its motion, such a robot can attain any state $(\theta,u,v)$ by only combining rotations and motion in the direction it is facing, i.e., by moving in only linear combinations of the $V_1$ and $V_2$ directions.
For example, imagine the robot has the state $(0,0,0)$, i.e., it is located at the origin and facing along the positive $u$-axis.
Suppose the robot rotates through an angle $\tau$, rolls forward at unit speed for time $\tau$, rotates through an angle $-\tau$, then rolls backward at unit speed for time $\tau$. 
This lands the robot at the point
\[
f(\tau) := \left(\tau\cos(\tau) - \tau,\tau\sin(\tau)\right)
\]
facing along the positive $u$-axis,
and
\[
\lim_{\tau \to 0} \tfrac{1}{\tau^2} f(\tau) = \lim_{\tau \to 0}\left(\frac{\cos(\tau) - 1}{\tau},\frac{\sin(\tau)}{\tau}\right) = (0,1).
\]
In other words, combining motions of order $\tau$ in the $V_1$ and $V_2$ directions followed by these motions in reverse results in movement in the $V_3$ direction on the order of $\tau^2$.
More precisely (see Definition~\ref{d-bracket}), on $\mathbb{R}^3$ we have
\begin{align*}
    [V_1,V_2](p) 
    := 
    DV_2(p)V_1(p) - DV_1(p)V_2(p) 
    = 
    \left( 
    \begin{array}{ccc}
       0  & 0 & 0 \\
       -\sin(\theta)  & 0 & 0 \\
       \cos(\theta) & 0 & 0
    \end{array}
    \right)
        \left(
    \begin{array}{c}
         1  \\
         0 \\
         0
    \end{array}
    \right)
    -
    0
    =
    V_3(p).
\end{align*}
That is, the prohibited vector field $V_3$ can be obtained by combining the vector fields $V_1$ and $V_2$ in this way.
In what follows, we will generalize this idea to collections of $m$ ``allowable'' directions in $\mathbb{R}^d$ with $d = m+\binom{m}{2}$.

\subsection{Lie brackets of vector fields}

Fix a positive integer $d$. For $i \in \{1,\dots,d\}$, we will use the notation $\frac{\partial}{\partial x^i}$ to denote the $i$th unit vector in $\mathbb{R}^d$, i.e., 
\[
\frac{\partial}{\partial x^i} = \mathbf{e}^i = (0,\dots,0,1,0,\dots,0).
\]
\begin{definition}[Smooth vector fields]
We define a {\em (smooth) vector field} to be a map $v:\mathbb{R}^d \to \mathbb{R}^d$
such that
\[
v(x) = \sum_{j=1}^d v^j(x) \frac{\partial}{\partial x^j}
=
(v^1(x),\dots,v^d(x))
\]
for all $x \in \mathbb{R}^d$
where the coefficients $v^j:\mathbb{R}^d \to \mathbb{R}$ are smooth functions.
\end{definition}

The following definition allows us to combine vector fields to form new ones as described in the setting of the two-wheeled robot allowing the robot to travel to any location by only moving in a prescribed subspace of directions.

\begin{definition}[Lie bracket]
\label{d-bracket}
Given two vector fields $v=\sum_{j=1}^d v^j \frac{\partial}{\partial x^j}$ and $w=\sum_{j=1}^d w^j \frac{\partial}{\partial x^j}$ on $\mathbb{R}^d$, we define the {\em Lie bracket} $[v,w]$ to be the vector field $[v,w]:\mathbb{R}^d \to \mathbb{R}^d$ defined as
\[
[v,w](x)
=
Dw(x)v(x) - Dv(x)w(x)
=
\sum_{i,j = 1}^d
\left(v^j(x) \frac{\partial w^i}{\partial x^j} (x) - w^j(x) \frac{\partial v^i}{\partial x^j}(x) \right)\frac{\partial}{\partial x^i}
\]
where $Dv$ and $Dw$ are the Jacobian matrices of $v$ and $w$ respectively.
\end{definition}
Note that, by definition, $[v,w]=-[w,v]$.
Heuristically, 
$Dw(x)v(x)$ is the infinitesimal translation of the vector $v(x)$ along the vector field $w$ at $x$.
If this infinitesimal translation is equal to the infinitesimal translation of  $w(x)$ along the vector field $v$ at $x$, then $[v,w]=0$, and we say that $v$ and $w$ {\em commute}.
Note, for example, that the coordinate vectors in $\mathbb{R}^d$ all commute with one another, i.e., $\left[ \frac{\partial}{\partial x^i},\frac{\partial}{\partial x^j}\right]=0$ for $i,j \in \{1,\dots,d\}$.

\subsection{Sub-Riemannian manifolds}
Fix $m < d$, and suppose $\mathcal{V} = \{v_1,\dots,v_m\}$ is a collection of smooth vector fields on $\mathbb{R}^d$.
Write
\[
\mathcal{V}_1(x) := \operatorname{span}\{
v_i(x):i \in \{1,\dots,m\}
\}
\]
and
\[
\mathcal{V}_2(x) := \operatorname{span}\{
[v_i,v_j](x):i,j \in \{1,\dots,m\}
\}
\]
for all $x \in \mathbb{R}^d$.

\begin{definition}[Sub-Riemannian manifold]
\label{d-SRM}
We say that $(\mathbb{R}^d,\mathcal{V})$ is an {\em equiregular sub-Riemannian manifold of step 2 with rank $m$}, or, more simply, a {\em sub-Riemannian manifold} if,
for every $x \in \mathbb{R}^d$, 
\begin{itemize}
    \item $\operatorname{dim}\mathcal{V}_1(x) = m$ for all $x \in \mathbb{R}^d$, i.e. $\{v_1(x),\dots,v_m(x)\}$ is a linearly independent set, and
    \item 
    $\mathcal{V}_1(x)
+
\mathcal{V}_2(x) = \mathbb{R}^d$.
\end{itemize}
\end{definition}

\subsubsection{Where this definition fits in the existing theory}
\label{sss-existingtheory}
In general, one defines a sub-Riemannian manifold to be a smooth, $d$-dimensional manifold $\mathcal{M}$ with a set of vector fields $\mathcal{V}$ on $\mathcal{M}$ that are {\em bracket generating} or {\em satisfy the H\"{o}rmander condition}, i.e., for any $x \in \mathcal{M}$, any other point can be written as a (finite) linear combination of vectors of the form
    \[
    [v_{i_1},v_{i_2}](x), \quad
    [v_{i_1},[v_{i_2},v_{i_3}]](x),
    \quad
    [v_{i_1},[v_{i_2},[v_{i_3}, v_{i_4}]]](x) \quad \cdots
    \qquad
    i_j \in \{1,\dots,m\}
    \]
    where $j \leq s$ for some $s \in \mathbb{N}$.
    The smallest such integer $s$ is called the {\em step} of $\mathcal{V}$ at $x$, and, if the step is the same at every point, we say that $\mathcal{V}$ is {\em equiregular} of step $s$.
    In other words, any vector in the tangent space of $\mathcal{M}$ at $x$ can be realized as a (finite) linear combination of vectors in $\mathcal{V}$ at $x$ and nested brackets of pairs of these vectors of length at most $s$.
    
    In our case, we consider $\mathcal{M} = \mathbb{R}^d$. Our second assumption in Definition~\ref{d-SRM} ensures that we need only consider a single bracket of vectors from $\mathcal{V}_1(x)$ to span all of $\mathbb{R}^d$ ({\em step 2}), and the first assumption ensures that $\mathcal{V}_1(x)$ is always $m$-dimensional ({\em rank $m$}) and that vectors from $\mathcal{V}_1(x)$ alone do not span $\mathbb{R}^d$ for any $x \in \mathbb{R}^d$ ({\em equiregular}).

        \subsection{Examples}
        Let us first observe some examples of vector field structures on $\mathbb{R}^d$. The first is a trivial example of a pair of vector fields on $\mathbb{R}^3$ that do not generate a sub-Riemannian structure. The other examples show choices of vector fields for which $\mathbb{R}^d$ is a sub-Riemannian manifold.

        \begin{figure}[h]
    \centering
    \includegraphics[width=0.31\linewidth]{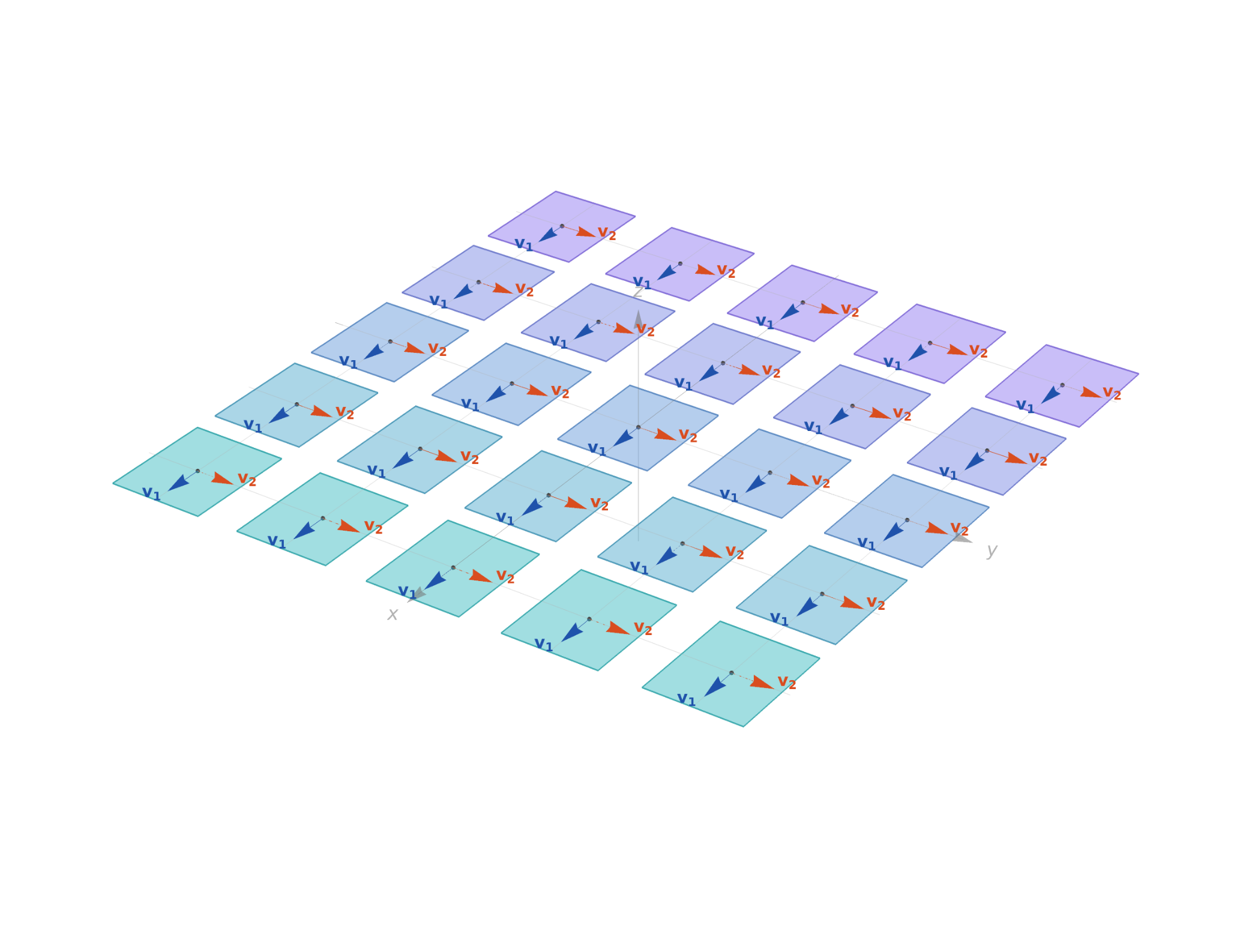}
    \includegraphics[width=0.31\linewidth]{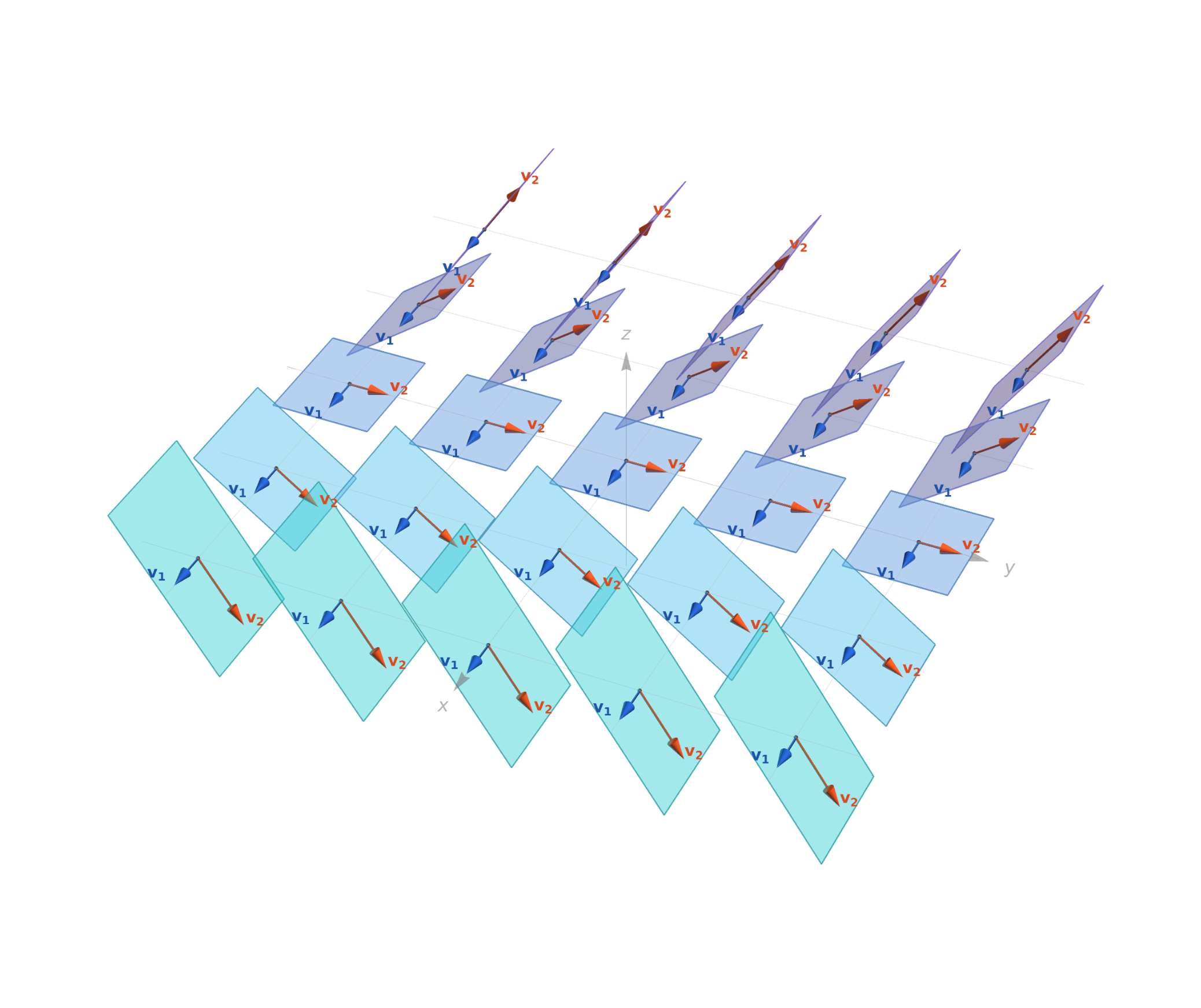}
    \includegraphics[width=0.31\linewidth]{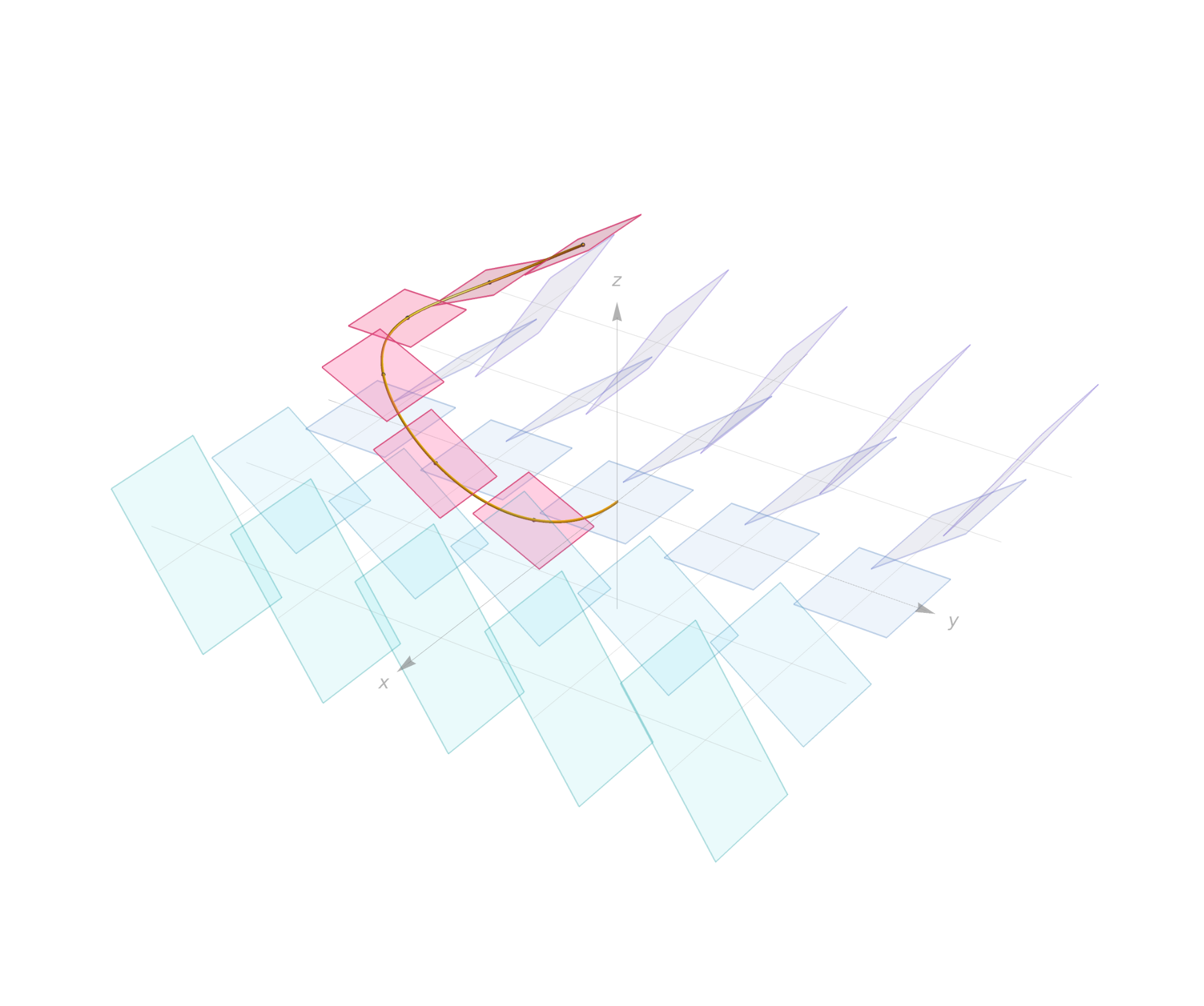}
    \caption{{\em Left}: Trivial horizontal planes $\mathcal{V}_1(x^1,x^2,0)$. \textit{Center}: Horizontal planes $\mathcal{V}_1(x,y,0)$ in the Heisenberg group $\mathbb{R}^3$. \textit{Right}: A horizontal curve in the Heisenberg group. Note that vertical lines are not horizontal.}
    \label{fig:planes}
\end{figure}

\subsubsection{Vector fields with rank 2 that are not bracket generating}
        Trivially, we may choose the coordinate vectors $v_1(x)=\frac{\partial}{\partial x^1}$ and $v_2(x)=\frac{\partial}{\partial x^2}$ on $\mathbb{R}^3$.
        See the leftmost image in Figure~\ref{fig:planes}.
        In this case, $\mathcal{V}_1(x) = \mathbb{R}^2 \times \{0 \}$ for all $x \in \mathbb{R}^3$, and so $\operatorname{dim}\mathcal{V}_1(x)=2$ for every $x \in \mathbb{R}^3$. However,
        \begin{align*}
    [v_1,v_2](x)
    =
    Dv_2(x)v_1(x)-Dv_1(x)v_2(x)
    &=
    0.
    \end{align*}
    Thus $\mathcal{V}_2(x)=\{0\}$, so
    $
    \mathcal{V}_1(x) + \mathcal{V}_2(x) = \mathcal{V}_1(x) 
    \neq \mathbb{R}^3$.

    \subsubsection{The Heisenberg group}
    \label{sss:heisenberg}
    Consider the example of the robot from the beginning of this section.
    Recall the vector fields on $\mathbb{R}^3$:
    \[
V_1(\theta,u,v) = \frac{\partial}{\partial \theta},
\qquad
V_2(\theta,u,v) = \cos(\theta) \frac{\partial}{\partial u} + \sin(\theta) \frac{\partial}{\partial v}.
\]
Let us approximate the behavior of the robot's motion near the angle $\theta =0$ so that $\cos \theta \approx 1$ and $\sin \theta \approx \theta$.
Making the change of coordinates $(x,y,z) \leftrightarrow (\theta,u,-v)$, we set
    \begin{equation}
    \label{e-HeisenbergVecs}
    v_1(x) = \frac{\partial}{\partial x}
    \qquad
    \text{and}
    \qquad v_2(x) = \frac{\partial}{\partial y}-
x \frac{\partial}{\partial z}. 
    \end{equation}
    These are the {\em Heisenberg vector fields} on $\mathbb{R}^3$.
    See Figure~\ref{fig:planes}.
    Moreover, writing $p=(x,y,z)$
    \begin{align*}
    [v_1,v_2](p)
    &=
    Dv_2(p)v_1(p)-Dv_1(p)v_2(p)\\
    &=
    \left( 
    \begin{array}{ccc}
       0  & 0 & 0 \\
       0  & 0 & 0 \\
       -1 & 0 & 0
    \end{array}
    \right)
    \left(
    \begin{array}{c}
         1  \\
         0 \\
         0
    \end{array}
    \right)
    -
        \left( 
    \begin{array}{ccc}
       0  & 0 & 0 \\
       0  & 0 & 0 \\
       0 & 0 & 0
    \end{array}
    \right)
    \left(
    \begin{array}{c}
         0  \\
         1 \\
         -x
    \end{array}
    \right)
    =     -\frac{\partial}{\partial z}
    \end{align*}
    and $[v_2,v_1](p) = \frac{\partial}{\partial z}$.
    In particular, $v_1$ and $v_2$ do not commute,
    and
    \begin{align*}
    \operatorname{span}\{v_1(p),v_2(p)\} + \operatorname{span}\{[v_1,v_2](p)\} 
    &= \operatorname{span}\left\{\tfrac{\partial}{\partial x} ,\tfrac{\partial}{\partial y} - x \tfrac{\partial}{\partial z}\right\} 
    +
    \operatorname{span}
    \left
    \{\tfrac{\partial}{\partial z}\right\}
    =
    \mathbb{R}^3.
    \end{align*}
    Therefore, the {\em Heisenberg group} $\mathbb H^1:=(\mathbb{R}^3, \{v_1,v_2\})$ is indeed a sub-Riemannian manifold (with step 2).

    \subsubsection{A direct sum of Heisenberg groups}
    \label{s-directSum}
    Fix a positive integer $n$, and set $m=2n$ and $d = 3n$.
    Express the coordinates of $\mathbb{R}^{3n}$ as
    \[
    p =
    (x^j,y^j,z^j)_{j=1}^n \in \mathbb{R}^{3n}.
    \]
    Write the $2n$ vector fields $\mathcal{V}=\{v_1^j,v_2^j\}_{j=1}^n$ as
        \begin{equation}
        \label{eqn:direct-sum-vfs}
    v_1^j(p) = \frac{\partial}{\partial x^{j}}
    \qquad
    \text{and}
    \qquad v_{2}^j(p) = \frac{\partial}{\partial y^{j}}-
x^{j} \frac{\partial}{\partial z^{j}}. 
    \end{equation}
When $n=1$, this is the standard Heisenberg group discussed above.
    The same calculation as above gives
        \begin{align*}
    [v_1^j,v_2^j](x)
    &=
    Dv_2^j(x)v_1^j(x)-Dv_1^j(x)v_2^j(x)
    =     -\frac{\partial}{\partial z^{j}},
    \end{align*}
    so $\mathcal{V}_2(x)= \operatorname{span}\{\frac{\partial}{\partial z^{j}} : j\in\{1,\dots,n\}\}$.
    Thus we can write
    \begin{align*}
        \mathbb{R}^d =
        \mathcal{V}_1(x) \oplus \mathcal{V}_2(x)
        =
        \bigoplus_{j=1}^n
        \left( \operatorname{span} \{v_1^j(x),v_2^j(x)\} \oplus
        \operatorname{span} \left\{\tfrac{\partial}{\partial z^{j}}\right\}
        \right).
    \end{align*}
    That is, $\mathbb{H}^n:=(\mathbb{R}^d,\mathcal{V})$ has the sub-Riemannian structure of the $n$-fold direct sum of the Heisenberg group with itself.

    \subsubsection{Engel algebra}

    The following vector fields have been studied extensively since the late 19th century. Recently, Morteza introduced the idea of flow matching using Continuous Normalizing Flow (CNF) on $\mathbb{R}^d$ with this sub-Riemannian structure (with slightly different notation) \cite{morteza2026parameter}.
    In the present paper, we will not consider this example (except when $d=3$) as the ODEs dictating geodesic flow these Engel groups are prohibitively computationally expensive.

    Define
    \[
        v_1(x) = \frac{\partial}{\partial x^1}
    \qquad
    \text{and}
    \qquad v_2(x) = \frac{\partial}{\partial x^2}-
x^{1} \frac{\partial}{\partial x^3} -
x^{3} \frac{\partial}{\partial x^4}
-
\cdots-
x^{d-1} \frac{\partial}{\partial x^d}. 
    \]
    Note that, when $d=3$, this is exactly the Heisenberg group. 
    When $d > 3$, the pair $(\mathbb{R}^d,\{v_1,v_2\})$ does not fit our definition of a sub-Riemannian manifold (Definition~\ref{d-SRM}) as it is no longer step 2. Indeed,
    \[
    \frac{\partial}{\partial x^{k+2}} = [\underbrace{v_2,[v_2,\cdots [v_2,[v_2}_{\text{$k$ copies}},v_1]] \cdots ]]
    \]
    for $k \in \{1,\dots,d-2\}$, so $(\mathbb{R}^d,\{v_1,v_2\})$ has step $d-1$.
    It does, however, fit the general definition of a sub-Riemannian manifold discussed in Section~\ref{sss-existingtheory}.

\subsection{Horizontal calculus}
    Fix $m <d$, let $\mathcal{V} = \{v_1,\dots,v_m\}$ be a collection of smooth vector fields on $\mathbb{R}^d$, and suppose $(\mathbb{R}^d,\mathcal{V})$ is a sub-Riemannian manifold.
    Recall the notation $\mathcal{V}_1(x) = \operatorname{span}\{
v_1(x),\dots,v_m(x)\}$ for $x \in \mathbb{R}^d$.
    \begin{definition}[Horizontal space]
    For each $x \in \mathbb{R}^d$, call
    \(
    \mathcal{V}_1(x) 
    \)
    the {\em horizontal space} at $x$.
    \end{definition}

        \begin{definition}[Horizontal frame]
    For each $x \in \mathbb{R}^d$, call the $d \times m$ matrix 
    $F(x)= [v_1(x) \cdots v_m(x)]$ 
    the {\em horizontal frame} at $x$.
    \end{definition}

    Observe that, for any vector $u = \sum_{i=1}^m u_iv_i(x) \in \mathcal{V}_1(x)$, if we write $u_H=(u_1,\dots,u_m) \in \mathbb{R}^m$, then $u = F(x)u_H$.
    
    \begin{definition}[Horizontal vector field]
    A vector field $w:\mathbb{R}^d \to \mathbb{R}^d$ is {\em horizontal at $x \in \mathbb{R}^d$} if 
    $
    w(x) \in \mathcal{V}_1(x)
    $.
    We say that $w$ is {\em horizontal} if it is horizontal at all $x \in \mathbb{R}^d$.
    \end{definition}
    In particular, a vector field $w$ is horizontal 
   if there are functions  $w_i: \mathbb{R}^d \to \mathbb{R}$ such that
    \[
    w(x) = \sum_{i=1}^m w_i(x) v_i(x).
    \]
    As a result, we will often use superscripts ($v^i$) to indicate Euclidean coordinates and subscripts ($v_i$) to indicate ``horizontal coordinates''.

    In the example of the two-wheeled robot, motion is only permitted in the directions lying in $\operatorname{span}\{V_1,V_2\}$, so the velocity field $v$ in the state space $\mathbb{R}^3$ corresponding to any allowable path the robot can travel must be a horizontal vector field. We call such paths the {\em horizontal curves}. See Figure~\ref{fig:planes} for an example of a horizontal curve in the Heisenberg group $\mathbb{R}^3$.
    
    \begin{definition}[Horizontal curve]
    An absolutely continuous curve $f:[0,1] \to \mathbb{R}^d$ is {\em horizontal} if 
    \begin{equation}
        \label{e-horizcurve}
        \dot{f}(t) \in \mathcal{V}_1(f(t)) \quad \text{for a.e. $t \in [0,1]$.}
    \end{equation}
    If we write $\dot f(t)=\sum_{i=1}^m a_i(t)v_i(f(t))$ a.e., then we call the map $a:[0,1] \to \mathbb{R}^m$ defined a.e. as $a(t)=(a_1(t),\dots,a_m(t))$ the {\bf horizontal control} of $f$.
    \end{definition}

    If $f$ is a horizontal curve with control $a$, then $\dot f(t) = F(f(t))a(t)$.

    Given a horizontal vector field $w$, write $w_H:=(w_1,\dots,w_m)$. Observe that 
    $w(x) = F(x)w_H(x)$ where $F(x)$ is the horizontal frame.
    In other words, for each $j \in \{1,\dots,d\}$, we have
        \begin{equation}
            \label{e-euctohorizvec}
            w^j = \sum_{i=1}^m w_i v_i^j.
        \end{equation}
    In particular, if $f:[0,1]\to \mathbb{R}^d$ is a horizontal curve and $a$ is its horizontal control, then 
    \begin{equation}
        \label{eqn:controlrelationship}
        \dot f(t) = V(f(t))a(t) \qquad \text{a.e. } t \in [0,1].
    \end{equation}

\subsection{Marginal horizontal vector fields}
\label{sec:marginalhorizontal}

Fix a set of $m$ linearly independent, bracket generating, left invariant vector fields $\{v_1,\dots,v_m\}$ on $\mathbb{R}^d$, and define the horizontal distribution as before: $\mathcal{V}_1(x) = \operatorname{span}\{v_1(x),\dots,v_m(x)\}$.

\begin{definition} 
We say that a map $u:[0,1] \times (\mathbb{R}^d)^3 \to \mathbb{R}^d$ is a {\em conditional, time dependent, horizontal vector field} if
$x \mapsto u(t,x,w,z)$ is a horizontal vector field on $\mathbb{R}^d$ for all fixed $t  \in [0,1]$ and $w,z \in \mathbb{R}^d$.
That is, $u_t(x|(w,z)) := u(t,x,w,z) \in \mathcal{V}_1(x)$ for all $x \in \mathbb{R}^d$.
\end{definition}

Given 
a state $x \in \mathbb{R}^d$, endpoints $w,z \in \mathbb{R}^d$, and a time $t \in [0,1]$, we write $\phi_t(x|(w,z)) := \phi_{x|(w,z)}(t)$ where $\phi_{x|(w,z)}:[0,1] \to \mathbb{R}^d$ is the curve satisfying
\begin{align*}
\dot{\phi}_{x|(w,z)}(t) &= 
u(t,\phi_{x|(w,z)}(t),(w,z))\\
 \phi_{x|(w,z)}(0) &= x
\end{align*}
That is, $\phi_t$ is the flow of the horizontal vector field $u_t(\cdot | (w,z))$, so, in particular, $t \mapsto \phi_t(x|(w,z))$ is a horizontal curve.

Fix a probability density function $p_0:\mathbb{R}^d \to (0,\infty)$, i.e., $\int_{\mathbb{R}^d}p_0 = 1$.
\begin{definition}
    We define a {\em conditional probability path} for the conditional, time dependent, horizontal vector field $u$ to be the function $p:[0,1]\times(\mathbb{R}^d)^3 \to (0,\infty)$ defined as
    follows. Fix $w,z \in \mathbb{R}^d$ and $t \in [0,1]$, and write $f(x):=\phi_t(x|(w,z))$ for all $x \in \mathbb{R}^d$ (so that $f$ is a diffeomorphism). Set
    \[
    p(t,x,w,z) := (f_{\#} p_0)(x) = p_0\left( f^{-1}(x) \right) \left|J(f^{-1})(x)\right|.
    \]
    We write $p_t(x|(w,z)) := p(t,x,w,z)$.
\end{definition}

Observe that $\int_{\mathbb{R}^d} p_t(x|(w,z)) \, dx = 1$
for all $t \in [0,1]$ and $w,z \in \mathbb{R}^d$.
We will consider $p_0$ a ``source'' density, and we 
fix a ``target'' probability density function $p_1:\mathbb{R}^d \to (0,\infty)$.
\begin{definition}
The {\em marginalized probability curve}
$\bar{p}:[0,1] \times \mathbb{R}^d \to (0,\infty)$ 
is defined for any 
$t \in [0,1]$ and $x\in\mathbb{R}^d$ as
\[
\bar{p}(t,x) := \int_{\mathbb{R}^d}\int_{\mathbb{R}^d} p_t(x|(w,z)) p_0(w)p_1(z) \, dw \, dz.
\]
Write $\bar{p}_t(x) := \bar{p}(t,x)$.
\end{definition}

Let us now define the marginalized vector field
$\bar{u}:[0,1] \times \mathbb{R}^d \to \mathbb{R}^d$ that $\bar{p}$ generates.
Fix $t \in [0,1]$ and endpoints $w,z \in \mathbb{R}^d$.
To construct $\bar{u}$, we average the $m$ horizontal coefficients $\mu_{t,i}$ of $u_t$ and reconstruct
the ambient marginal field using the frame $F$. Note that the endpoint
integration still ranges over $\mathbb R^d \times \mathbb{R}^d$. That is, we only reduce the number of independent control coefficients, not the size of the integration domain.
As we will see in Proposition~\ref{prop:SameAsChen}, the resulting field agrees with ambient componentwise
averaging as in \cite{chen2024flow}.
Recall that, since $u_t(\cdot|(w,z))$ is a horizontal vector field for any $t \in [0,1]$ and $w,z \in \mathbb{R}^d$,
there are real valued functions $\mu_{t,i}(\cdot|(w,z)):\mathbb{R}^d \to \mathbb{R}$ 
for $i \in \{1,\dots,m\}$
such that 
\[
u_t(x|(w,z)) = \sum_{i=1}^m \mu_{t,i}(x|(w,z))v_i(x).
\]
\begin{definition}[Marginalized Vector Field]
For each $i \in \{1,\dots,m\}$, define the function $\bar{\mu}_i:[0,1]\times \mathbb{R}^d \to \mathbb{R}$ as
\[
\bar{\mu}_i(t,x) := \int_{\mathbb{R}^d}\int_{\mathbb{R}^d} \mu_{t,i}(x|(w,z)) \frac{p_t(x|(w,z))p_0(w)p_1(z)}{\bar{p}_t(x)} \, dw \, dz,
\]
and write $\bar{\mu}(t,x):=(\bar{\mu}_1(t,x),\dots,\bar{\mu}_m(t,x))$.
Define the \textbf{marginalized vector field}
$\bar{u}:[0,1] \times \mathbb{R}^d \to \mathbb{R}^d$
as
\[
\bar{u}(t,x) 
:= 
F(x)\bar \mu(t,x)
=
\sum_{i=1}^m \bar{\mu}_{i}(t,x)v_i(x),
\]
and write $\bar{u}_t(x):=\bar{u}(t,x)$.
\end{definition}

The next proposition shows that averaging the control coefficients
and then applying the horizontal frame gives the same ambient
marginal field as componentwise averaging. In particular,
marginalization preserves horizontality.
\begin{proposition}
\label{prop:SameAsChen}
For all $t \in [0,1]$ and $x \in \mathbb{R}^d$, we have
    \[
    \bar{u}_t(x) = \sum_{j=1}^d \left(\int_{\mathbb{R}^d} \int_{\mathbb{R}^d} u_t^j(x|(w,z)) \frac{p_t(x|(w,z))p_0(w)p_1(z)}{\bar{p}_t(x)} \, dw \, dz\right) \frac{\partial}{\partial x^j}.
    \]
\end{proposition}
\begin{proof}
    As before, write $v_i = \sum_{j=1}^d v_i^j \frac{\partial}{\partial x^j}$.
    Recall from \ref{e-euctohorizvec} that $u_t^j(x|(w,z)) = \sum_{i=1}^m \mu_{t,i}(x|(w,z))v_i^j(x)$ for all $x \in \mathbb{R}^d$.
    Therefore,
    \begin{align*}
        \bar{u}_t(x) &= 
        \sum_{i=1}^m \bar{\mu}_i(t,x)v_i(x)\\
        &=
        \sum_{i=1}^m \left(\int_{\mathbb{R}^d} \int_{\mathbb{R}^d} \mu_{t,i}(x|(w,z)) \frac{p_t(x|(w,z))p_0(w)p_1(z)}{\bar{p}_t(x)} \, dw \, dz \right) \sum_{j=1}^d v_i^j(x)\frac{\partial}{\partial x^j}\\
        &=
        \sum_{j=1}^d \left(\int_{\mathbb{R}^d} \int_{\mathbb{R}^d} \left[ \sum_{i=1}^m \mu_{t,i}(x|(w,z)) v_i^j(x) \right] \frac{p_t(x|(w,z))p_0(w)p_1(z)}{\bar{p}_t(x)} \, dw \, dz \right) \frac{\partial}{\partial x^j}\\
        &=
        \sum_{j=1}^d \left(\int_{\mathbb{R}^d} \int_{\mathbb{R}^d} u_t^j(x|(w,z)) \frac{p_t(x|(w,z))p_0(w)p_1(z)}{\bar{p}_t(x)} \, dw \, dz\right) \frac{\partial}{\partial x^j}.
    \end{align*}
\end{proof}

We conclude this section by showing that the mass conservation law holds exactly as in \cite{chen2024flow}.
Note that $\bar{p}_t\bar{u}_t = \sum_{i=1}^m \bar{p}_t \bar{\mu}_i(t,\cdot)v_i$ is a horizontal vector field, so 
$\operatorname{div}_H(\bar{p}_t\bar{u}_t)$ is well defined.

\begin{theorem}[Mass Conservation]
If $\operatorname{div}(v_i)=0$ for $i = 1,\dots,m$, then
    \[
    \frac{\partial}{\partial t} \bar{p}_t = -\operatorname{div}(\bar{p}_t\bar{u}_t)= -\operatorname{div}_H(\bar{p}_t\bar{u}_t).
    \]
\end{theorem}
\begin{proof}
    The first equality follows exactly as in \cite{chen2024flow} by the previous proposition,
    and the second equality follows from the facts that $p_tu_t$ is horizontal and that our assumption implies $\operatorname{div} = \operatorname{div}_H$.
\end{proof}

\section{Free Carnot groups of step 2}
\label{sec:freecarnotsection}

In this section, we introduce and explore the particular sub-Riemannian structure on $\mathbb{R}^d$ that will be used for path generation.

Fix a positive integer $m$, and set $d = m + \binom{m}{2}$.
Since $d-m = \binom{m}{2}$, we can identify each $\ell \in \{m+1,\dots,d\}$
with a pair $(i,j)$ such that $i < j$.
Thus we may write $x \in \mathbb{R}^d$ using coordinates
\[
x=(x^1,x^2,\dots,x^m,x^{(1,2)},x^{(1,3)},\dots,x^{(1,m)},x^{(2,3)},x^{(2,4)},\dots,x^{(2,m)},\dots,x^{(m-1,m)}).
\]
Given $x \in \mathbb{R}^d$, define for $j = 1,\dots,m$
\begin{equation}
\label{e:step2vfs}
v_j(x) = \frac{\partial}{\partial x^j} - \sum_{i=1}^{j-1} x^{i} \frac{\partial}{\partial x^{(i,j)}},
\end{equation}
and write $\mathcal{V}= \{v_1,\dots,v_m\}$. Note that $v_i(x)$ and $v_j(x)$ are linearly independent for any $i,j \in \{1,\dots m\}$,
so $\mathcal{V}_1(x) = \operatorname{span}\{v_1(x),\dots,v_m(x)\}$ has dimension $m$.
Moreover, for each $i,j \in \{1,\dots,m\}$ with $i<j$, one may check that
\begin{align*}
[v_i,v_j] 
= -\frac{\partial}{\partial x^{(i,j)}}.
\end{align*}
Thus, for all $x \in \mathbb{R}^d$,
\begin{align*}
\mathcal{V}_2(x) &= 
\operatorname{span}\left\{[v_i,v_j]
:
i,j\in\{1,\dots,m\}, \, i<j
\right\} 
= \{0\} \times \mathbb{R}^{d-m}.
\end{align*}
In other words, $\mathbb{F}_{m,2} :=(\mathbb{R}^d,\mathcal{V})$ is a sub-Riemannian manifold as defined in Definition~\ref{d-SRM}.
$\mathbb{F}_{m,2}$ is sometimes called the {\em free Carnot group of step 2 with $m$ generators}. Observe also that $\mathbb{F}_{2,2}$ is exactly the Heisenberg group $\mathbb{H}^1$.

\subsection{Horizontal coordinates}
According to the definition of the horizontal vector fields \eqref{e:step2vfs}, we have 
$
v_i^i(x) = 1
$,
and
$
v_i^j(x) = 0
$ for all other coordinates $j \in \{1,\dots,m\}$ in the first layer.
Thus, for any horizontal vector $u =\sum_{i=1}^m u_iv_i(x)= F(x)u_H$, equation \eqref{e-euctohorizvec} gives
\[
u^j = \sum_{i=1}^m u_iv_i^j(x) = u_j
\]
for all $j \in \{1,\dots,m\}$ and all $x \in \mathbb{R}^d$.
In other words, the first $m$ Euclidean coordinates of $u$ are exactly the $m$ horizontal coordinates of $u$.
Thus, if $\pi_1:\mathbb{R}^d \to \mathbb{R}^m$ is the projection onto the first $m$ coordinates of $\mathbb{R}^d$, we have
\[
\pi_1(u) = (u_1,\dots,u_m) = u_H.
\]
Note that this is independent of the choice of $x$.
This gives us the following natural definition.
\begin{definition}[Horizontal inner product and length]
For every $x \in \mathbb{R}^d$ and any $u_1,u_2 \in \mathcal{V}_1(x) \subset \mathbb{R}^d$, define the {\em horizontal inner product} $\langle u_1 , u_2 \rangle_H : = \langle \pi_1(u_1),\pi_1(u_2)\rangle_{\mathbb{R}^m}$ and {\em horizontal length} $\Vert u_1 \Vert_H := |\pi_1(u_1)|_{\mathbb{R}^m}$.
\end{definition}

\subsection{Horizontal curves and geodesics}
The following is the benefit of the sub-Riemannian setting:
\begin{proposition}
\label{p:horizstep2Carnot}
Suppose $\mathbb{F}_{m,2} = (\mathbb{R}^d,\mathcal{V})$ is the sub-Riemannian manifold defined by \eqref{e:step2vfs}.
A smooth curve $\gamma:[0,1] \to \mathbb{R}^d$ is horizontal if and only if
\begin{align*}
    \dot{\gamma}^{(i,j)} 
    =
    - \dot{\gamma}^j\gamma^{i}.
\end{align*}
\end{proposition}

We may interpret this result as follows: given a horizontal curve $\gamma=(\gamma_H,\gamma_V):[0,1] \to \mathbb{R}^m \times \mathbb{R}^{d-m}$, the coordinates of $\gamma_V$ are determined {\em uniquely} (up to initial value $\gamma_V(0)$) by $\gamma_H$ and $\dot{\gamma}_H$. Indeed, the Fundamental Theorem of Calculus gives for all $t \in [0,1]$
    \begin{equation}
    \label{e-lift}     
\gamma^{(i,j)}(t) 
= \gamma^{(i,j)}(0)-\int_0^t \dot{\gamma}^j\gamma^{i}.
\end{equation}

    \begin{proof}[Proof of Proposition~\ref{p:horizstep2Carnot}]
    Since $\gamma$ is horizontal, we have
    for some real valued functions $g_j$ for $j \in \{1,\dots,m\}$ that
        \begin{align*}
    \dot{\gamma} 
    = 
    \sum_{j=1}^m g_j v_j(\gamma)
    =
    \sum_{j=1}^m g_j \left(\frac{\partial}{\partial x^j} - \sum_{i=1}^{j-1} \gamma^{i} \frac{\partial}{\partial x^{(i,j)}}\right)
    =
    \sum_{j=1}^m g_j \frac{\partial}{\partial x^j} - \sum_{j=1}^m \sum_{i=1}^{j-1} g_j\gamma^{i} \frac{\partial}{\partial x^{(i,j)}}.
\end{align*}
Since also 
\[
\dot{\gamma} = \sum_{j=1}^m \dot{\gamma}^j \frac{\partial}{\partial x^j} - \sum_{j=1}^m \sum_{i=1}^{j-1} \dot{\gamma}^{(i,j)} \frac{\partial}{\partial x^{(i,j)}},
\]
it must be the case that $g_j = \dot{\gamma}^j$, from which the result follows.
    \end{proof}

\begin{definition}
[Horizontal length]
Given a horizontal curve $\gamma:[0,1] \to \mathbb{R}^d$, define
\[
\ell(\gamma) 
:= 
\int_0^1 \Vert \dot{\gamma}\Vert_H =
\int_0^1 \sqrt{  (\dot{\gamma}^1)^2+ \cdots + (\dot{\gamma}^m)^2}.
\]
We call $\ell(\gamma)$ the {\bf horizontal length} of $\gamma$ and $\Vert \dot{\gamma}\Vert_H$ its {\bf horizontal speed}.
\end{definition}

\subsection{Group structure}

\begin{definition}[Group operation]    
Given a point $x \in \mathbb{R}^d$, define $L_x:\mathbb{R}^d \to \mathbb{R}^d$ as
\begin{equation}
    \label{e-groupop}
    \begin{aligned}
(L_x( y))^j &= x^j+y^j,
&& 1\le j\le m,\\
(L_x(y))^{(i,j)}
&= x^{(i,j)}+y^{(i,j)}-x^i y^j,
&& 1\le i<j\le m.
\end{aligned}
\end{equation}
\end{definition}
One may show that the binary operation $x * y := L_x(y)$ is actually a group operation on $\mathbb{R}^d$, i.e., it is closed, associative, has an identity element (namely, the origin $0 \in \mathbb{R}^d$), and every point $x \in \mathbb{R}^d$ has an inverse $x^{-1} \in \mathbb{R}^d$. In particular,
\begin{equation}
\label{e-groupinv}
    \begin{aligned}
(x^{-1})^j &= -x^j,
&& 1\le j\le m,\\
(x^{-1})^{(i,j)}
&= -x^{(i,j)}-x^i x^j,
&& 1\le i<j\le m.
\end{aligned}
\end{equation}
Note however, that this group operation is {\em not} abelian. That is, in general, we have
\[
x*y \neq y*x,
\qquad
\text{i.e., } L_x(y) \neq L_y(x).
\]
Moreover, as $L_x$ acts linearly in the first $m$ coordinates, we may conclude that
\[
\ell(x*\gamma) = \ell(L_x \circ \gamma) = \ell(\gamma)
\]
for any horizontal curve $\gamma$ and any $x \in \mathbb{R}^d$.

\subsection{Regression of horizontal controls }

Our goal is to predict a random horizontal control $c \in \mathbb{R}^m$ given a time $t$ and state $x$ using a predictor $b(t,x) \in \mathbb{R}^m$. The error in this prediction can be measured by using the norm $| \cdot |_{\mathbb{R}^m}$ within the horizontal space $\mathcal{V}_1$ spanned by the orthonormal vectors $\{v_1,\dots,v_m\}$ or by passing the controls to the ambient space $\mathbb{R}^d$ using the frame $F(x)$ and measuring with the norm $\vert \cdot \vert_{\mathbb{R}^d}$.
That is, we define the loss functions as
\[
 L_H(b)=\mathbb E|b(t,x)-c|_{\mathbb{R}^m}^2,\qquad
 L(b)=\mathbb E\left|F(x)b(t,x)-F(x)c\right|_{\mathbb{R}^d}
 ^2
\]
for any measurable predictor $b:\mathbb{R} \times \mathbb{R}^d \to \mathbb{R}^m$.
According to the next result, these loss functions are both minimized by the same predictor. This minimizing predictor is the conditional expectation $\bar c(t,x)=\mathbb E[C\mid T=t,X=x]$.

\begin{proposition}[Horizontal and ambient loss functions]
\label{thm:loss-comparison}
Let $(T,X,C)\in \mathbb{R} \times \mathbb{R}^d \times \mathbb{R}^m$ be a random time, state, and horizontal
control such that
$\mathbb E| C |_{\mathbb{R}^m}^2$
and
$\mathbb E|F(X)C|_{\mathbb{R}^d}^2$ are finite.
Then
\begin{equation}
\label{eqn:lossesaremnimnized}
L_H(\bar c)=\min L_H(b)
\qquad
\text{ and }
\qquad
L(\bar c) = \min L(b)
\end{equation}
where the minimum is taken over all measurable predictors $b:\mathbb{R} \times \mathbb{R}^d \to \mathbb{R}^m$.
Moreover, any other predictor $\hat b$ satisfying 
\eqref{eqn:lossesaremnimnized}
must have $\hat b=\bar c$ almost everywhere with respect to the distribution of \((T,X)\).
\end{proposition}
While $L_H$ and $L$ are both minimized by $\bar c$, their values $L_H(\bar c)$ and $L(\bar c)$ may be different.

\begin{proof}
For a fixed $(t,x)$, write $\mu_{(t,x)}$ to denote the distribution on the controls $C$. 
Writing
$b-C=(b-\bar c)+(\bar c-C)$ then gives
\begin{align}
\label{eqn:innerproduct-integral-loss}
    \int|b-C|_{\mathbb{R}^m}^2 \, d\mu_{(t,x)}(C) = \int \left(|b-\bar c|_{\mathbb{R}^m}^2 + 2 \langle b- \bar c, \bar c - C\rangle_{\mathbb{R}^m} + |\bar c - C|_{\mathbb{R}^m}^2\right) \, d\mu_{(t,x)}(C).
\end{align}
Since $(t,x)$ is fixed, the definition of $\bar c$ gives
\[
\int \langle b- \bar c, \bar c - C\rangle_{\mathbb{R}^m} \, d\mu_{(t,x)}(C)
=
\left \langle b-c,\bar c - \int C \, d\mu_{(t,x)}(C) \right \rangle_{\mathbb{R}^m}
=0.
\]
Therefore, among all possible choices of predictors $b$, \eqref{eqn:innerproduct-integral-loss} is minimized by the choice $b= \bar c$. Integrating this equation over all pairs $(t,x)$ gives the second equality in \eqref{eqn:lossesaremnimnized}.

The first equality and the a.e.-uniqueness of the minimizer follow from the fact that the vectors that make up the frame $F(x)$ are orthonormal. Indeed, for a fixed $x \in \mathbb{R}^d$, 
\[
\langle z_1,z_2 \rangle_{F(x)} := z_1^T F^T(x)F(x)z_2,
\qquad
z_1,z_2 \in \mathbb{R}^m
\]
is an inner product on $\mathbb{R}^m$ satisfying 
$\langle F(x)z_1,F(x)z_2 \rangle_{\mathbb{R}^d} = \langle z_1,z_2 \rangle_{F(x)}$,
so we can repeat the above argument replacing the inner product $\langle\cdot,\cdot\rangle_{\mathbb{R}^m}$
with $\langle\cdot,\cdot\rangle_{F(x)}$.
\end{proof}

\section{Candidate curves in the free Carnot groups}
\label{sec:candidatecurves}

We describe two constructions of horizontal paths between
given endpoints in the setting of the free Carnot groups $\mathbb{F}_{m,2} = (\mathbb{R}^d,\mathcal{V})$ just discussed. 
These are the paths that will be used in training the network for SRCFM.
The first computes sub-Riemannian paths
using the geodesic ODE in~\eqref{e-newODEforfree}. Note that a solution to this equation need not be a geodesic curve in general.
The second is an explicit piecewise curve consisting of
a horizontal straight segment followed by a
loop constructed from a spectral decomposition of a skew-symmetric matrix. Both constructions supply horizontal training paths, and the numerical validity of their endpoints will be assessed separately.

\begin{figure}[h]
    \centering
    \includegraphics[width=0.9\linewidth]{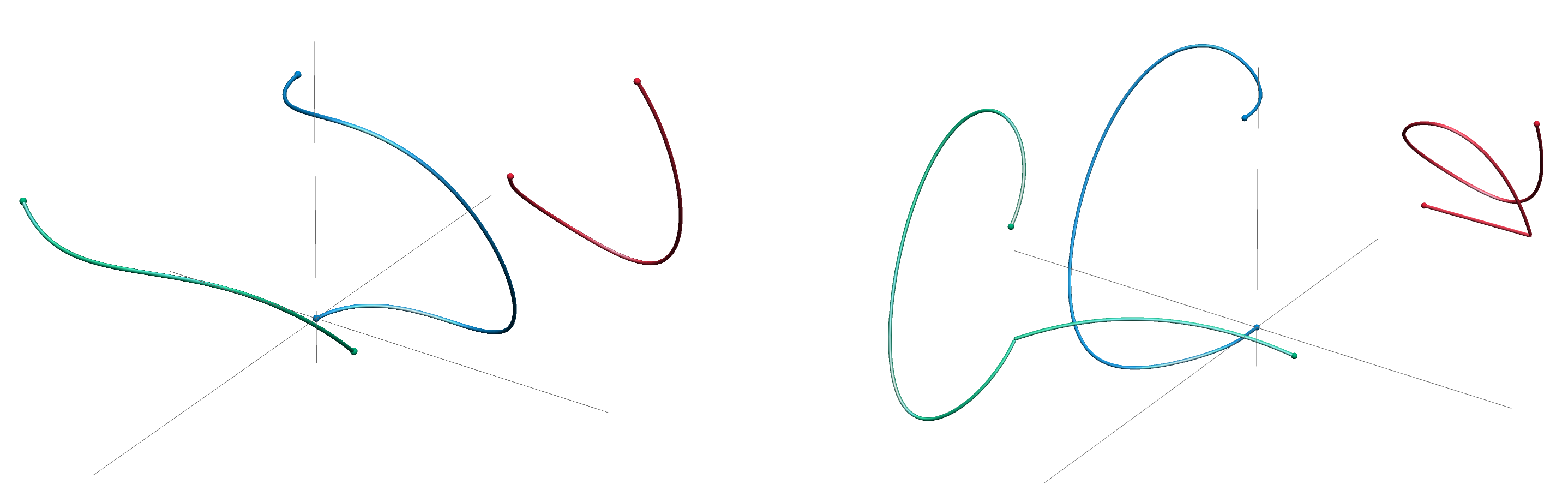}
    \caption{Training paths in the Heisenberg group $\mathbb{R}^3$.  The left-hand figure shows paths generated using the ODE in \eqref{e-newODEforfree}, while the explicit piecewise curves from \eqref{eq:single-excursion-path} are on the right.
    }
    \label{fig:placeholder}
\end{figure}

\subsection{General sub-Riemannian paths: endpoint equations}
\label{sec:geodesic}

\begin{definition}
[Length-minimizing geodesic]
Fix $x_0,x_1 \in \mathbb{R}^d$.
A horizontal curve $\gamma:[0,1] \to \mathbb{R}^d$ is a \textbf{length minimizing geodesic} (or simply a {\bf geodesic}) from $x_0$ to $x_1$ if its horizontal length
$
\ell(\gamma)
$
is minimal among all horizontal curves from $x_0$ to $x_1$.
\end{definition}
Of course, there is no initial guarantee that a length minimizing curve exists, but, under the assumptions in Definition~\ref{d-SRM}, such a curve is guaranteed to exist (see, for example, \cite[Theorem 3.43]{agrachev2020}).

Suppose $\gamma:[0,1] \to \mathbb{R}^d$ is a geodesic curve from $0$ to $q \in \mathbb{R}^d$. 
Write $\gamma=(\gamma_H,\gamma_V):[0,1] \to \mathbb{R}^m \times \mathbb{R}^{d-m}$.
It follows from a Lagrange multipliers argument that the horizontal part $\gamma_H$ of $\gamma$ must satisfy
\begin{equation*}
\ddot{\gamma}_H = -A \dot{\gamma}_H
\end{equation*}
for some skew symmetric $m \times m$ matrix $A$.
(See, for example, the discussion preceding Theorem~2 in \cite{brockett1980}.)
Integrating this ODE, we may equivalently write that $\gamma$ must satisfy 
\begin{equation}
    \label{e-newODEforfree}
\dot{\gamma}_H = v_0-A \gamma_H
\end{equation}
where $v_0 \in \mathbb{R}^m$ is the initial velocity of $\gamma_H$.
It follows that $M_A(t)v_0$ solves \eqref{e-newODEforfree} where $M_A(t) :=\int_0^t e^{-sA} \, ds$.
This gives us the following.
See Equation (13) in \cite{rizziserres2017}.
\begin{theorem}
\label{t-geodformula}
    If $\gamma:[0,1] \to \mathbb{R}^d$ is a geodesic with $\gamma(0)=0$ and $\dot\gamma(0)=v_0$,
    then there is a matrix
    $A \in \mathfrak{so}(m)$ satisfying $\gamma_H(t)=M_A(t)v_0$.
    In particular, if $\gamma(1) = q$, then $A$ satisfies the IVP
    \begin{align}
    \label{e-Xequation}
    (M_A(1)v_0)^j = q^j,
    \qquad
    \left(-\int_0^1 
    \left[M_A(t)v_0\right] 
\left[
\dot{M}_A(t)v_0
\right]^T  dt \right)^{(i,j)}
=
q^{(i,j)},
\qquad
1\leq i < j \leq m.
\end{align}
    \end{theorem}
    Equation \eqref{e-Xequation} follows immediately from \eqref{e-lift}.
    Note that a horizontal curve satisfying this IVP is not guaranteed to be a geodesic, only that geodesics must satisfy the IVP.

\subsubsection{Geodesics in the Heisenberg group and its direct sums}
\label{sec:specialgeodesiccase}

Let us investigate the special case of the Heisenberg group described in Section~\ref{sss:heisenberg}, i.e., $d= 3$ and $m=2$.
The following discussion is well known in the literature (see, e.g., \cite{hajlaszZimmerman2015}).
Write $0\neq q=(x,y,z) \in \mathbb{R}^3$, so $q_H = (x,y)$ and $q_V=z$. 

Suppose first that $q_H\neq 0$ and $q_V = z \neq 0$.
Fix a non-zero skew symmetric matrix $A \in \mathfrak{so}(2)$.
It must be true that
$A =\left(\begin{smallmatrix}
0 & \lambda\\
-\lambda & 0
\end{smallmatrix}\right)$
for some real parameter $\lambda$,
so, as long as $A \neq 0$
\[
M_A(t) 
= 
\frac{1}{\lambda}
\begin{pmatrix}
\sin(t\lambda) & \cos(t\lambda)-1\\
1-\cos(t\lambda) & \sin(t\lambda)
\end{pmatrix}.
\]

Suppose $\gamma=(f,g,h):[0,1] \to \mathbb{R}^3$ is a geodesic from 0 to $q = (x,y,z) \in \mathbb{R}^3$ with $(x,y) \neq (0.0)$.  Theorem~\ref{t-geodformula} guarantees the existence of a matrix $A=\left(\begin{smallmatrix}
    0 & \lambda \\ -\lambda &0
\end{smallmatrix}\right)$ satisfying \eqref{e-Xequation}.
Suppose $\lambda \in (-2\pi, 2 \pi)$.
If $\lambda \neq 0$, then $M_A(t)$ is invertible, so $v_0 = M_A(1)^{-1}(x,y)$.
Substituting this into \eqref{e-Xequation} gives
\[
-z = 
\left(\int_0^1 
    \left[M_A(t)M_A(1)^{-1}(x,y)\right] 
\left[
\dot{M}_A(t)M_A(1)^{-1}(x,y)
\right]^T  dt \right)^{(1,2)}
=
\frac{xy}{2} +
\frac{(x^2+y^2)(\lambda-\sin\lambda)}
{4(1-\cos\lambda)}.
\]

If we write 
$
\kappa(t) := \frac{\sin t - t}
{2(1-\cos t)}
$,
then the task of determining the matrix $A$
when $\lambda \neq 0$
reduces to solving
\begin{equation}
\label{e:flambda}    
\kappa(\lambda) = \frac{2z + xy}{x^2+y^2}
\end{equation}
for $\lambda$.
Note that, if $\lambda = 0$, then $A=0$, so $M_A(t)=t I_{2 \times 2}$.
Finally, if $(x,y)=(0,0)$, set $\lambda = 2 \pi \operatorname{sgn}(z)$.

In any case, we have found a matrix $A=\left(\begin{smallmatrix}
    0 & \lambda \\ -\lambda &0
\end{smallmatrix}\right)$ 
such that $(f(t),g(t)) = M_A(t)v_0$, and $v_0 = M_A(1)^{-1}(x,y)$.
Defining $h(t)=-\int_0^t \dot g f$ gives the explicit formula for the geodesic $\gamma=(f,g,h)$.

In the setting of the direct sum of Heisenberg groups described in
Section~\ref{s-directSum}
(where
$m=2n$ and $d=3n$ for a given positive integer $n$),
we may apply the above argument in each triplet $(x^j,y^j,z^j)$ to obtain a curve $\gamma=(f^j,g^j,h^j)_{j=1}^n:[0,1] \to \mathbb{R}^{3n}$.
That is, for $j \in \{1,\dots,n\}$, solve $\kappa(\lambda_j)=\frac{2z^j+x^jy^j}{(x^j)^2+(y^j)^2}$ for $\lambda_i$, and obtain $n$ separate curves $(f^j,g^j,h^j)$.

Note for future reference that the group law and inverses are computed independently in each triplet. That is, if we write $p_k=(x_k^j,y_k^j,z_k^j)_{j=1}^n$ for $k=1,2,3$ and set $p_3 = p_1 * p_2$, then 
\begin{equation}
    \label{e-groupopHeisenberg}
    \begin{aligned}
x_3^j &= x_1^j + x_2^j\\
y_3^j &= y_1^j + y_2^j \qquad \qquad \qquad 1 \leq j \leq n, \\
z_3^j &= z_1^j + z_2^j - x_1^jy_2^j \\
\end{aligned}
\end{equation}
and
\begin{equation}
\label{e-groupinvHeisenberg}
    \begin{aligned}
    \left((x^j,y^j,z^j)_{j=1}^n\right)^{-1}
    =
    (-x^j,-y^j,-z^j-x^jy^j)_{j=1}^n.
\end{aligned}
\end{equation}

\subsection{Piecewise sub-Riemannian paths: explicit formulas}

Recall that we can write any skew symmetric matrix $A \in \mathfrak{so}(m)$ as
\begin{equation}
\label{eq:single-excursion-decomposition}
 A=\sum_{k=1}^{r}\lambda_k
 (a_k b_k^T-b_k a_k^T),\qquad
 \lambda_1\geq\cdots\geq\lambda_r>0,\qquad
 r\leq\lfloor m/2\rfloor
\end{equation}
where $a_k \pm ib_k$ are the complex eigenvectors of $A$ with eigenvalues $\pm i\lambda_k$.

The following theorem (specifically, the formula \eqref{eq:single-excursion-path}) gives us instructions to build a horizontal curve connecting any two points in $\mathbb{R}^d$. Notably, this curve will have constant horizontal speed almost everywhere.

\begin{theorem}
\label{thm:single-excursion-connector}
Fix $q=(q_H,q_V)\in \mathbb{R}^m \times \mathbb{R}^{d-m}$ with $q \neq 0$, and define the skew-symmetric matrix $A$ by
\begin{equation}
\label{eq:single-excursion-area}
 A_{ij}=q^{(i,j)}+\tfrac12q^i q^j\quad(i<j),\qquad
 A_{ji}=-A_{ij},\qquad A_{ii}=0.
\end{equation}
Using the decomposition \eqref{eq:single-excursion-decomposition}, we define for $s \in [0,1]$
\begin{equation}
\label{eq:single-excursion-eta}
 \eta(s)=\sum_{k=1}^{r}\sqrt{\frac{\lambda_k}{\pi k}}
 \bigl[(\cos(2\pi k s)-1)a_k-\sin(2\pi k s)b_k\bigr],
 \qquad
 K(s)=\int_0^s\eta\dot\eta^T.
\end{equation}
Define $\tau=|q_H|/L$ 
and $s(t)=(t-\tau)/(1-\tau)$ for $t \in (\tau,1]$ 
where 
\begin{equation}
\label{eq:single-excursion-length}
 L=|q_H| + \sqrt{4\pi\sum_{k=1}^{r}k\lambda_k}.
\end{equation}
Define the curve
$\gamma=(h,v):[0,1] \to \mathbb{R}^m \times \mathbb{R}^{d-m}$ as
\begin{equation}
\label{eq:single-excursion-path}
\begin{aligned}
 h(t)&=
 \begin{cases}
  (t/\tau)q_H,&0\leq t\leq\tau,\quad \tau>0,\\
  q_H+\eta(s(t)),&\tau\leq t\leq1,\quad \tau<1,
 \end{cases}\\[2mm]
 v^{(i,j)}(t)&=
 \begin{cases}
  -\dfrac{t^2}{2\tau^2}q^i q^j,
       &0\leq t\leq\tau,\quad \tau>0,\\
  -\tfrac12q^i q^j-q^i\eta^j (s(t))-[K(s(t))]^{(i,j)},
       &\tau\leq t\leq1,\quad \tau<1,
 \end{cases}
 \qquad i<j.
\end{aligned}
\end{equation}
Then $\gamma$ is horizontal
with
$\gamma(0)=0$ and $\gamma(1)=q$, and
$|\dot h(t)| = L$ for all $t \neq \tau$.

For arbitrary distinct points $x_0,x_1\in \mathbb{R}^d$, set
$\Gamma=L_{x_0}\circ\gamma$. Then $\Gamma$ is a horizontal
curve from $x_0$ to $x_1$ with constant horizontal speed.
\end{theorem}

\begin{proof}
We first note that $\gamma$ is horizontal since
\[
 \dot v^{(i,j)}(t)=-\frac{t}{\tau^2}q^i q^j
 =-h^i(t)\dot h^j(t)
\]
for all $t \in [0,\tau)$, and
\[
 \dot v^{(i,j)}(t)
 =-q^i \tfrac{d}{dt} \eta^j(s(t))-\tfrac{d}{dt}K(s(t))^{(i,j)}
 =-(q^i+ \eta^i(s(t))) \tfrac{d}{dt}\eta^j(s(t))
 =-h^i(t)\dot h^j(t)
\]
for $t \in (\tau,1]$.
Also, it is clear that $\gamma(0)=0$.

It remains to prove that $\gamma(1)=q$ and $|\dot h(t)| = L$ whenever $t \neq \tau$.
Since $\eta(1) = 0$, it follows that $h(1) = q_H$.
Thus, to prove that $\gamma(1)=q$, we must show that 
$v^{(i,j)}(1) = q_V^{(i,j)}$.
Observe that
\begin{equation}
    \label{e-goalofheight}
    v^{(i,j)}(1) + \tfrac12 q^iq^j = - [K(1)]^{(i,j)}.
\end{equation}
According to \eqref{eq:single-excursion-area} and the definition of $K$, our goal is therefore to prove that $K(1) = \int_0^1 \eta \dot \eta^T = - A$.

To do this, write for each $k \in \{1,\dots,r\}$
\[
 \eta_k(s)=\sqrt{\frac{\lambda_k}{\pi k}}
 \bigl[(\cos(2\pi k s)-1)a_k
       -\sin(2\pi k s)b_k\bigr]
\]
so that 
\[
\int_0^1 \eta \dot \eta^T =
\int_0^1 \left(\sum_{k=1}^r \eta_k\right) \left( \sum_{\ell=1}^r \dot\eta_\ell^T\right)
=
\int_0^1 \sum_{k=1}^r \sum_{\ell=1}^r \eta_k\dot\eta_\ell^T.
\]
For any $k,\ell \in \{1,\dots,r\}$,
\[
\begin{aligned}
 2\int_0^1\cos(2\pi\ell s)\cos(2\pi k s)\,ds
 &=\delta_{\ell k}, \qquad 
 2\int_0^1\sin(2\pi\ell s)\sin(2\pi k s)\,ds
 =\delta_{\ell k},\\
 \int_0^1\cos(2\pi\ell s)\sin(2\pi k s)\,ds&= \int_0^1 \cos(2\pi\ell s) \, ds = \int_0^1 \sin(2\pi\ell s)\,ds=0
\end{aligned}
\]
where $\delta_{\ell k}=1$ when $\ell=k$ and $\delta_{\ell k}=0$ otherwise.
Therefore, $\int_0^1 \eta_\ell\dot\eta_k^T = 0$ for all $\ell \neq k$. 

When $\ell=k$, the coefficients of
$a_k a_k^T$ and $b_k b_k^T$ are constant multiples of $\cos(2\pi k s)\sin(2\pi k s)$ and $\sin(2 \pi k s)$.
Therefore,
\begin{align*}
    K(1) 
    = \int_0^1 \eta_k \dot\eta_k^T 
    &=\sum_{k=1}^r \int_0^1 \eta_k \dot\eta_k^T \\
    &= 
    \sum_{k=1}^r 2 \lambda_k\left[\int_0^1(\cos(2\pi ks) - \cos^2(2\pi ks))a_kb_k^T \, ds
    +
    \int_0^1\sin^2(2\pi ks)b_ka_k^T \, ds
    \right]
    \\
    &=
    -\sum_{k=1}^r\lambda_k (a_k b_k^T-b_ka_k^T) = -A.
\end{align*}
Combining this with \eqref{e-goalofheight} gives us $\gamma(1) = q$ as desired.

Finally, 
for all $t \in (0,\tau)$, we have $|\dot h(t)| = |q_H|/\tau = L$, and for $t \in (\tau,1)$, $|\dot h(t)| = |\frac{d}{dt}\eta(s(t))| = \frac{|\dot \eta(s(t))|}{1-\tau}$.
Since the vectors $\{a_1,b_1,\dots,a_r,b_r\}$ are orthonormal, we have
for $s \in (0,1)$
\[
 |\dot \eta(s)|^2
 =
 \sum_{k=1}^r|\dot \eta_k(s)|^2
 =\sum_{k=1}^{r}4 \pi k \lambda_k
   \bigl(\sin^2(2 \pi k s)+\cos^2(2 \pi k s)\bigr)
 =4\pi\sum_{k=1}^{r}k \lambda_k.
\]
The definitions of $\tau$ and $L$ then give $\frac{|\dot\eta|}{1-\tau} = L$ as desired.

For arbitrary distinct endpoints $x_0,x_1 \in \mathbb{R}^d$, 
define $q \in \mathbb{R}^d$ as
\begin{equation*}
\label{eq:single-excursion-relative}
 q^j=x_1^j-x_0^j,\qquad
 q^{(i,j)}=x_1^{(i,j)}-x_0^{(i,j)}+x_0^i q^j.
\end{equation*}
The group law gives 
$q=x_0^{-1}*x_1$.
Suppose $\gamma=(h,v)$ is the curve from 0 to $q$ defined above, and set $\Gamma = L_{x_0} \circ \gamma$.
Then
\[
 \Gamma^j=x_0^j+h^j,\qquad
 \Gamma^{(i,j)}=x_0^{(i,j)}+v^{ij}-x_0^i h^j.
\]
Consequently, $\dot\Gamma^j=\dot h^j$ and
$\dot\Gamma^{(i,j)}=-(x_0^i+h^i)\dot h^j=-\Gamma^i \dot{\Gamma}^j$,
so $\Gamma$ is also horizontal. The fact that $\Gamma$ connects $x_0$ to $x_1$ follows from $\gamma(0)=0$ and $\gamma(1)=q=x_0^{-1}*x_1$.
\end{proof}

\section{Methodology and Algorithms}

We use the two-endpoint formulation of CFM with independent
sampling~\cite[Section~3.2.2]{tong2024improving}.
Let $\mu$ and $\nu$ be the source and target distributions on $\mathbb{R}^d$.
We draw $w\sim\mu$ and $z\sim\nu$ independently, so the endpoint pair
$(w,z)$ has joint law $\pi=\mu\otimes\nu$. The marginal probability path
and velocity field are obtained by marginalizing over both endpoints,
as in Section~\ref{sec:marginalhorizontal}.

For each pair $(w,z)$, our goal is to choose a curve $\gamma_{w,z}:[0,1]\to \mathbb{R}^d$ connecting them. In the standard Euclidean case, we typically use the straight line path
\begin{equation}
\label{eqn:straightline}    
 \gamma_{w,z}(t)=(1-t)w+tz,\qquad
 a_{w,z}(t)=z-w.
\end{equation}
We then train a neural network by sampling endpoints $(w,z)\sim\pi$ and time $t\sim\mathrm{Unif}[0,1]$, and we regress the predicted velocity at $(t,\gamma_{w,z}(t))$ against the goal of $a_{w,z}$.
This is the zero-noise case of the two-endpoint construction from \cite{tong2024improving}.

SRCFM uses the same method described above, but it replaces the straight line curve \eqref{eqn:straightline} and its velocity curve with horizontal curves $\gamma_{w,z}(t) \in \mathbb{R}^d$ and their {\bf horizontal controls} $a_{w,z}(t) \in \mathbb{R}^m$.
In what follows, we will describe three different routines 
for selecting $\gamma_{w,z}$ and $a_{w,z}$ using the discussion from Section~\ref{sec:candidatecurves}.
Though we use different routines to select the curves,
the
training and generation algorithm will be the same; varying
the paths changes the regression targets, not the training procedure.

\subsection{Computing paths between two endpoints}
\label{sec:geodesic-target-algorithms}

A path routine takes endpoints $w,z \in \mathbb{R}^d$ and query times $t \in [0,1]$, and it returns a curve $\gamma_{w,z}(t) \in \mathbb{R}^d$ and its horizontal control $a(t) \in \mathbb{R}^m$.  To do so, use the appropriate choice of group law and inverse to write $q=w^{-1}*z$, and use the discussion from previous sections (as outlined below) to construct a path $\gamma_q$ from $0$ to $q$ with control $a_q$. Finally, translate to obtain the path $\gamma_{w,z}=w * \gamma_q$ with control $a_{w,z}$.  According to the definition of the group operation \eqref{e-groupop}, this translation leaves the horizontal control unchanged, so, in fact, $a_{w,z}=a_q$.  Let us now explain how to obtain $\gamma_q$ and $a_q$ for a given $q$.

\subsubsection{Heisenberg group and direct sums.}
\label{sss:Heisenbergpaths}
Fix the Heisenberg geometry $(\mathbb{H}^1)^n$ on $\mathbb{R}^{3n}$, and use the group law and inverse in \eqref{e-groupopHeisenberg} and \eqref{e-groupinvHeisenberg}.
Write $q=(x^j,y^j,z^j)_{j=1}^n$, and construct $\gamma_q$ by applying the construction described in Section~\ref{sec:specialgeodesiccase}
to each triple of coordinates at the same time parameter $t$. 

\subsubsection{General sub-Riemannian paths.}
\label{sss:SRpaths}
Fix the free step 2 geometry $\mathbb{F}_{m,2}=(\mathbb{R}^d,\mathcal{V})$. Thus $d = m + \binom{m}{2}$, and we use the group law and inverse in \eqref{e-groupop} and \eqref{e-groupinv}.  Write $M_A(t)=\int_0^t e^{-sA}\,ds$ for any $m \times m$ matrix $A$.  Keeping $v_0 \in \mathbb{R}^m$ and the skew-symmetric matrix $A \in \mathfrak{so}(m)$ as unknowns, solve \eqref{e-Xequation} by damped Newton or nonlinear least squares to obtain a curve $h(t)=M_A(t)v_0$.  Define $v(t)$ using \eqref{e-lift} with $v(0)=0$, and return the curve $\gamma_q=(h,v)$ with and $a_q(t)=e^{-tA}v_0$. Note as before that this curve is not guaranteed to be a globally length minimizing geodesic but only that every such geodesic must satisfy \eqref{e-Xequation}.

\subsubsection{Piecewise sub-Riemannian paths.}
\label{sss:piecewisepaths}
We again use the free step 2 geometry $\mathbb{F}_{m,2}$.
Evaluate Theorem~\ref{thm:single-excursion-connector} directly, and construct the curve $\gamma_q=(h,v)$ using
\eqref{eq:single-excursion-area}--\eqref{eq:single-excursion-path}.
If $A=0$, use $\tau=1$, and ignore the second case in \eqref{eq:single-excursion-path}. If $q_H=0$, use $\tau = 0$ and ignore the first case. 
Set $a_q=\dot{h}$. This is defined for $t \neq \tau$, and we can set $a_q(\tau)=q_H$. Uniform time sampling will hit this value with probability zero.

\begin{algorithm}[htbp]
\caption{Path evaluation using the preceding constructions}
\label{alg:srcfm-endpoint-path}
\begin{algorithmic}[1]
\Require Endpoints $w,z$, query times, path choice and associated geometry (including any retained pair set) \ref{sss:Heisenbergpaths}--\ref{sss:piecewisepaths}. Endpoint tolerance $\varepsilon$.
\State Calculate $q=w^{-1}*z$ with the chosen group law
\If{$q=0$}
 \State \Return the constant path $\gamma_q = w$ and $a_q=0$
\EndIf
\If{the choice is \ref{sss:Heisenbergpaths}}
 \State Apply Section~\ref{sec:specialgeodesiccase} with the limiting cases described above
 \State \Return $\gamma_q=(f^j,g^j,h^j)_{j=1}^n$ and $a_q=(\dot f^j, \dot g^j)_{j=1}^n$
\ElsIf{the choice is \ref{sss:SRpaths}}
 \State Solve \eqref{e-Xequation} for $(A,v_0)$ within the prescribed solver budget
 \State \Return $\gamma_q=(M_A(t)v_0,v)$ and $a_q=e^{-tA}v_0$
 \State If no endpoint solution is accepted, save diagnostics and \Return failure
\Else
 \State Use Theorem~\ref{thm:single-excursion-connector} 
 \State \Return $\gamma_q=(h,v)$ and $a_q=\dot h$
\EndIf
\State Evaluate $\gamma_q$ and $a_q$ at the query times
\State Set $\gamma_{w,z} = w * \gamma_q$ using the chosen group law 
\State Set $a_{w,z}=a_q$ 
\State Check reconstructed and independently integrated endpoints against $w,z$
\State If either check exceeds $\varepsilon$, save diagnostics and \Return failure
\State \Return $\gamma_{w,z}(t),a_{w,z}(t)$ at the requested times
\end{algorithmic}
\end{algorithm}

\subsection{One training and generation algorithm for all path choices}
\label{sec:srcfm-method}
\label{sec:two-endpoint-cfm}

For each pair $(w,z)$, choose a horizontal path
$\gamma_{w,z}$ with its prescribed endpoints and control $a_{w,z}$, so
$\dot\gamma_{w,z}(t)=F(\gamma_{w,z}(t))a_{w,z}(t)$ almost everywhere
where $F(x)$ is the frame of the chosen geometry. (See \eqref{eqn:controlrelationship}.)
For $(W,Z)\sim\pi$, Proposition~\ref{thm:loss-comparison} tells us that the regression target for square-integrable controls should be
\[
 \bar a(t,x)=\mathbb E[a_{W,Z}(t)\mid \gamma_{W,Z}(t)=x],\qquad
 V(t,x)=F(x)\bar a(t,x).
\]
Thus the average $\bar a$ is taken over paths $\gamma_{W,Z}$ passing through the same state at the same time for varying endpoints $W,Z$. 
For a fixed $t$ and $x$, one can first average over the starting point $W$ conditional on the current state $\gamma_{W,Z}(t)=x$ and the terminal point $Z$, then average over terminal points $Z$
conditional on $\gamma_{W,Z}(t)=x$.

The neural network observes pairs $(t,x)$ where $x$ lies in the full $d$-dimensional
state and predicts only the $m$-dimensional control $a_\theta(t,x)\in\mathbb R^m$. These predicted controls then form a time dependent horizontal vector field $V_\theta$ and generate flows $\Phi^\theta_t$ as follows:
\[
 V_\theta(t,x)=F(x)a_\theta(t,x),\quad
  \dot\Phi^\theta_t(w)=V_\theta(t,\Phi^\theta_t(w)),\quad
 \Phi^\theta_0(w)=w.
\]
Therefore, the endpoint map $T_\theta(w):=\Phi^\theta_1(w)$ is the terminal ODE map wherever the solution exists.
While $T_\theta$ is learned from endpoint pairs, it is evaluated only at a source point. 
Given $(w,z)\sim\pi$ and $t\sim\mathrm{Unif}[0,1]$, the control regression objective is
\begin{equation}
\label{eq:srcfm-objective}
 \mathcal L_H(\theta)=
 \mathbb E
 |a_\theta(t,\gamma_{w,z}(t))-a_{w,z}(t)|_{\mathbb{R}^m}^2.
\end{equation}
Alternatively, one may consider loss in the ambient space
\begin{equation}
\label{eq:srcfm-ambient-objective}
 \mathcal L(\theta)=
 \mathbb E|F(\gamma_{w,z}(t))
 [a_\theta(t,\gamma_{w,z}(t))-a_{w,z}(t)]|_{\mathbb{R}^d}^2.
\end{equation}
Proposition~\ref{thm:loss-comparison} guarantees that these expectations are both minimized when $a_\theta$ is exactly $\bar a$, but they have potentially different behavior in the setting of finite-network approximation. The chosen loss and any normalization must therefore be specified in practice.

\begin{algorithm}[htbp]
\caption{SRCFM training and generation with a chosen path routine}
\label{alg:srcfm-training}
\begin{algorithmic}[1]
\Require Laws $\mu,\nu$, coupling $\pi$, frame $F$, compatible routine $\operatorname{Path}_F$
\Require Network $a_\theta$, batch size $B$, times per pair $K$, optimizer, ODE solver for generation
\For{each training update}
 \State Draw fresh endpoint pairs $(w_i,z_i)$ according to $\pi$, $i=1,\ldots,B$
 \State Draw independent $t_{ij}\sim\mathrm{Unif}[0,1]$, $j=1,\ldots,K$
 \For{$i=1,\ldots,B$}
  \State Construct and validate one path for $(w_i,z_i)$ at all requested times
  \State If construction fails, save diagnostics and stop; do not replace the pair silently
  \State Obtain $(g_{ij},a_{ij})=\operatorname{Path}_F(w_i,z_i;t_{ij})$
 \EndFor
 \State Update $\theta$ using $(BK)^{-1}\sum_{i,j}\|a_\theta(t_{ij},g_{ij})-a_{ij}\|^2$
\EndFor
\State Draw fresh $w\sim\mu$ and integrate $\dot g=F(g)a_\theta(t,g)$, $g(0)=w$
\State \Return $T_\theta(w)=g(1)$ and, if required, its trajectory; report any integration failure
\end{algorithmic}
\end{algorithm}

With $K=1$ this is the usual minibatch estimator; $K>1$ reuses each
endpoint solve at several independent times. Given a pair of endpoints, all chosen times use
the same path parameters including the chosen endpoint solution or
decomposition of the skew-symmetric matrix. Alternate path constructions or
endpoint-dependent solvers can be used in this algorithm if declared and measurable.
However, the path family must remain fixed during the computation of each training target.

\begin{remark}[Conditional and generated paths]
\label{rem:srcfm-not-geodesic}
Whine the averaging of the paths at a fixed point preserves horizontality, it need not
preserve the geodesic equations even when the conditional paths themselves are
geodesics. This mirrors the setting of Euclidean CFM in which the straight segment conditional paths may give rise to curved marginal paths.
Similarly, the piecewise training paths need not result in piecewise generated
marginalized paths. 
\end{remark}

\section{Numerical results}

\label{sec:numerical-results}

We first study transport between distributions on Heisenberg groups $\mathbb{H}^1$, their direct sums $(\mathbb{H}^1)^n$,
and on the free step-two group structure $\mathbb{F}_{m,2}$ from Section~\ref{sec:freecarnotsection}. 
In each example, we generate
samples by integrating the learned velocity field and compare their
endpoints with the target distribution. We also measure how closely
the numerical trajectories satisfy the horizontal equations. We then
compare paths between fixed endpoints, test the piecewise construction
with 100 horizontal controls, and use it to train SRCFM on resized
MNIST images.

For the step-two coordinates of Section~3, we write $h_j=x^j$ and
$v_{ij}=x^{(i,j)}$. The horizontal controls $a_j=\dot h_j$ determine the
vertical velocities by $\dot v_{ij}=-h_i a_j$, as in \eqref{e-lift}.

\subsection{Evaluation measures}
\label{sec:num-protocol}
\label{sec:num-metrics}

For two sets of $M$ points, generated samples $\widehat x_i$ and target
samples $y_i$, we measure the empirical Wasserstein distance
\[
 \widehat W_2=
 \left(\frac1M\min_{\sigma\in S_M}
 \sum_{i=1}^M\|\widehat x_i-y_{\sigma(i)}\|_2^2\right)^{1/2}.
\]
Here $S_M$ is the set of permutations of $M$ indices. This measure uses
the best one-to-one matching in Euclidean distance over all state
coordinates. For the Heisenberg tests, we also report biased squared
maximum mean discrepancy (MMD). For two samples of size $M$, it compares
kernel similarities through
\[
 \widehat{\mathrm{MMD}}_{\mathrm b}^{\,2}
 =\frac1{M^2}\sum_{i,j=1}^M
 \left[k(\widehat x_i,\widehat x_j)+k(y_i,y_j)
       -2k(\widehat x_i,y_j)\right].
\]
We use
$k(u,v)=\sum_{\ell\in\{0.5,1,2\}}\exp(-\|u-v\|_2^2/(2\ell^2))$.

To measure horizontality, let $g_i^r$ be the state of trajectory $i$ at
saved time $t_r$. We estimate the velocity on each interval by
$Dg_i^r=(g_i^{r+1}-g_i^r)/(t_{r+1}-t_r)$ and evaluate the coefficients of
the horizontal equations at $\overline g_i^r=(g_i^{r+1}+g_i^r)/2$.
For $N$ trajectories on $J$ intervals, we average the absolute errors
in the vertical equations. In the direct sum of $n$ Heisenberg groups, the horizontal residual is
\[
 \mathcal R_{\mathbb H}
 =\frac1{NJn}\sum_{i=1}^N\sum_{r=0}^{J-1}\sum_{k=1}^n
 \left|Dz_{i,k}^r-(-\overline x_{i,k}^rDy_{i,k}^r)\right|.
\]
The residual includes the effects of time integration and finite
differences. It uses every vertical equation, even when a figure
shows only a few coordinates.

\subsection{Single and product Heisenberg experiments}
\label{sec:num-heisenberg}

In the first example, we transport a Gaussian source to a noisy helix
on one Heisenberg group and then to a product of ten such helices.
This lets us compare the distributions generated by CFM and SRCFM
and measure how closely their trajectories satisfy the prescribed
horizontal equations. We take $n=1$ or $n=10$ factors as in
Section~\ref{s-directSum}. In each factor, 
the horizontality condition is 
$\dot z=-x\dot y$. The full state thus has $3n$
coordinates and the horizontal control has $2n$ components.

The source law is $\mu=\mathcal N(0,0.49I_{3n})$. We obtain the target
by sampling each factor independently according to
\[
 \begin{gathered}
 \vartheta\sim\mathrm{Unif}[0,2\pi],\qquad R=1.60+0.18\epsilon_R,\\
 (x,y,z)=(R\cos\vartheta+0.08\epsilon_x,\,
 R\sin\vartheta+0.08\epsilon_y,\,
 0.50\vartheta+0.20\epsilon_z),
 \end{gathered}
\]
where the four $\epsilon$ variables are independent standard Gaussians
and independent of $\vartheta$. As the angle $\vartheta$ increases, the
mean position turns in the $(x,y)$ plane and rises in $z$. The Gaussian
perturbations spread the samples around this curve and give the target
a density in the ambient coordinates.

At each update, we draw fresh source and target batches. We sort each
batch by the mean vertical coordinate $n^{-1}\sum_{k=1}^n z_k$ and
pair samples by their positions in the sorted lists. We draw training
times independently and uniformly on $[0,1]$. CFM uses straight
segments between paired endpoints, and SRCFM solves the Heisenberg
curvature equation separately in each factor by 30 bisection iterations
with curvature magnitude bounded by $2\pi-10^{-6}$. In the curvature
calculation, we floor the squared horizontal displacement at $10^{-10}$
and set the curvature to zero when the displacement squared is below
this threshold or when the magnitude of the ratio in the curvature
equation is below $10^{-7}$. We use the straight-path formula when
the curvature magnitude is below $10^{-6}$ and floor its associated
determinant in the path calculation at $10^{-12}$. We report these safeguards because they affect near-vertical
training paths, whose endpoint errors were not recorded.

Both networks take time and the full state as inputs. Each has three
hidden layers of width 256 with GELU activations and a linear output
layer. The CFM network predicts all $3n$ velocity components. The
SRCFM network predicts the $2n$ horizontal controls and obtains the
vertical velocities from the horizontal equations. We minimize squared
error averaged over samples and output coordinates. Table~\ref{tab:num-setup}
gives the resulting network sizes.

\begin{table}[th!]
\centering\small
\caption{State dimensions, numbers of horizontal controls, and network
parameter counts in the Heisenberg experiments.}
\label{tab:num-setup}
\begin{tabular}{lrrrr}
\toprule
Geometry & State & Controls & CFM parameters & SRCFM parameters\\
\midrule
$\mathbb H^1$ & 3 & 2 & 133,635 & 133,378\\
$(\mathbb H^1)^{10}$ & 30 & 20 & 147,486 & 144,916\\
\bottomrule
\end{tabular}
\end{table}

We train on a GPU for 6,000 updates using Adam with batch size 256,
constant learning rate $10^{-3}$, coefficients $(0.9,0.999)$, and
$\epsilon=10^{-8}$. Weight decay and gradient clipping are disabled.
For each method we evaluate the final network from five training runs.
In each comparison, the two methods use the same 600 source samples
and 600 independent target samples. We integrate the learned fields
using 100 fourth-order Runge--Kutta (RK4) steps on $[0,1]$. We compare the
first 300 generated and target samples in the original coordinates
and measure the horizontal residual over all 600 paths.

\begin{figure}[th!]
\centering
\includegraphics[width=0.72\textwidth]{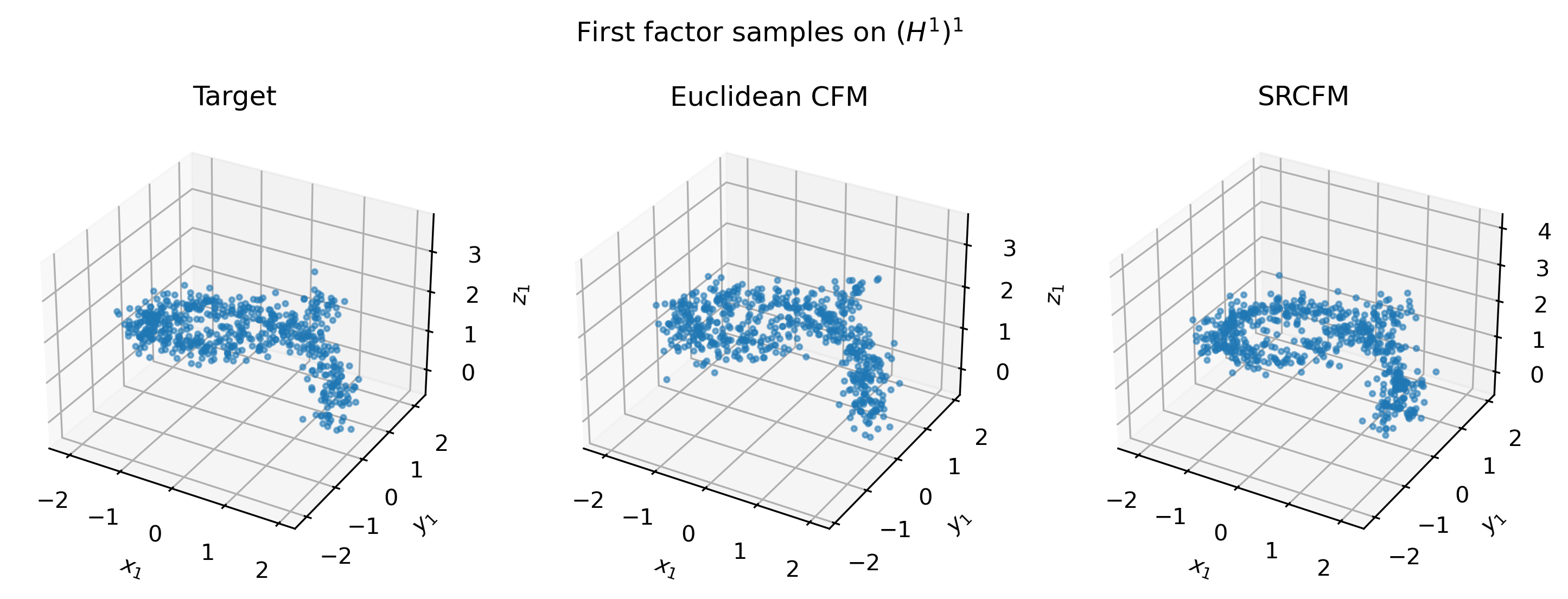}\\[4pt]
\includegraphics[width=0.72\textwidth]{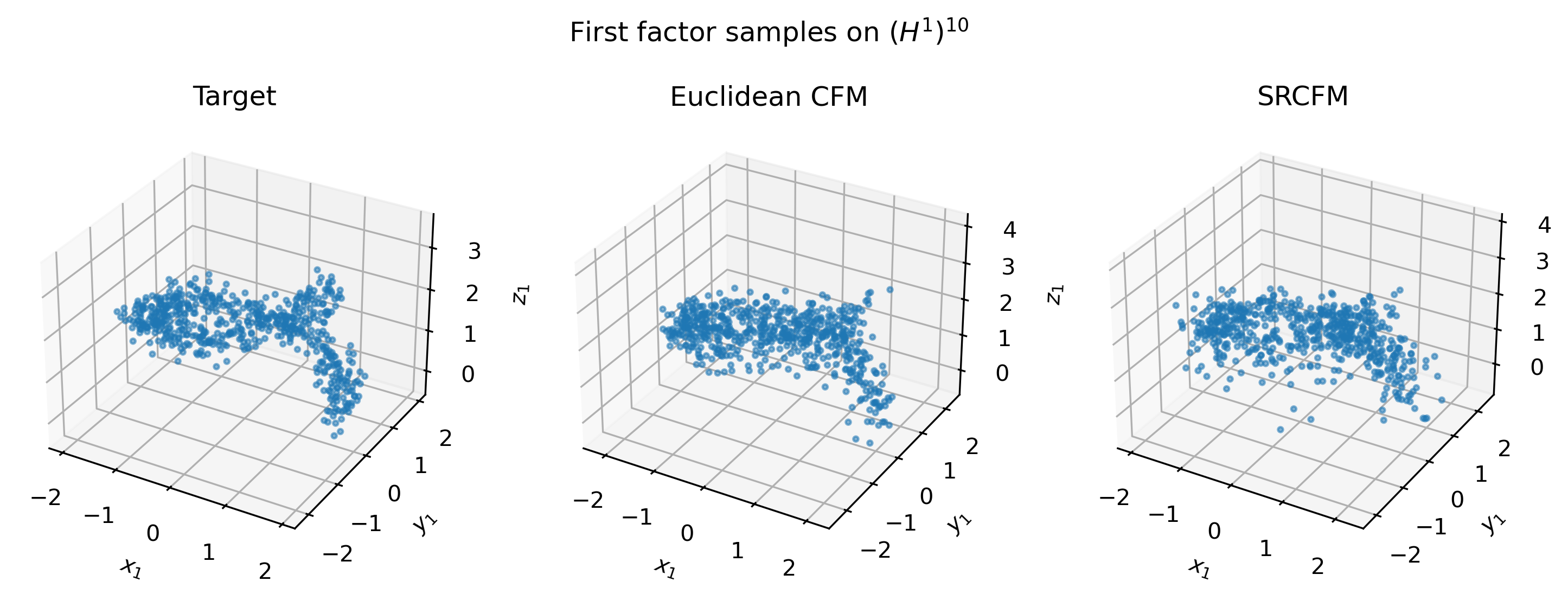}
\caption{Target samples and generated endpoints from one run per method. The top
row shows $\mathbb H^1$. The bottom row shows the first factor of
$(\mathbb H^1)^{10}$. Each panel contains 600 points and uses
independently scaled axes.}
\label{fig:num-heisenberg-samples}
\end{figure}

Figure~\ref{fig:num-heisenberg-samples} shows the target samples on the
left and the endpoints generated by CFM and SRCFM in the center and
right columns. In the top row, the two generated clouds reproduce the
main curved band of the target, including the bend near its lower end. In the bottom row, we see the same overall bend in
the first factor of the ten-group product, with more points scattered
away from the target band. This row shows a three-dimensional view of
the generated samples, so we use the full thirty-dimensional states
to assess agreement with the target law.

Table~\ref{tab:num-results} shows slightly lower mean Wasserstein
errors for CFM and much smaller horizontal residuals for SRCFM.
In the product case, squared MMD adds little to the comparison
because both values are close to the contribution from comparing
each point with itself.

\begin{table}[tb]
\centering\small
\setlength{\tabcolsep}{4pt}
\caption{Heisenberg results, reported as mean $\pm$ sample standard
deviation over five runs. Distribution errors use 300 generated and
300 target samples. The horizontal residual $\mathcal R_{\mathbb H}$
is evaluated against the Heisenberg horizontal equations for all 600
paths generated by each method.}
\label{tab:num-results}
\begin{tabular}{lccc}
\toprule
Method & $\widehat W_2$ &
$\widehat{\mathrm{MMD}}_{\mathrm b}^{\,2}$ & $\mathcal R_{\mathbb H}$\\
\midrule
\multicolumn{4}{l}{\textbf{Single Heisenberg group} ($\mathbb H^1$, $d=3$)}\\
CFM & $0.3617\pm0.0966$ &
$0.01331\pm0.00393$ & $1.5373\pm0.0505$\\
SRCFM & $0.3778\pm0.0731$ &
$0.01555\pm0.00630$ & $(2.605\pm0.175)\,10^{-4}$\\
\addlinespace
\multicolumn{4}{l}{\textbf{Product of ten Heisenberg groups} ($(\mathbb H^1)^{10}$, $d=30$)}\\
CFM & $5.4273\pm0.0438$ &
$0.020036\pm0.000058$ & $1.5887\pm0.0251$\\
SRCFM & $5.5480\pm0.0290$ &
$0.020068\pm0.000018$ & $(1.944\pm0.049)\,10^{-4}$\\
\bottomrule
\end{tabular}
\end{table}

\subsection{\texorpdfstring{Gaussian transport on the free group $\mathbb{F}_{4,2}=\mathbb{R}^{10}$}{Gaussian transport on the free group F4,2=R10}}
\label{sec:num-geodformula-verification}
\label{sec:num-single-excursion}
\label{sec:num-gaussian-general}

In the second example, we use the paths associated with
Theorem~\ref{t-geodformula} as training targets on $\mathbb{F}_{4,2} = \mathbb{R}^{10}$. We choose
\[
 \mu=\mathcal N(0,0.5^2I_{10}),\qquad
 \nu=\mathcal N(c_\star,0.3^2I_{10}),\qquad
 c_\star=0.75(e_{h_1}+e_{v_{12}}),
\]
where $e_{h_i}$ and $e_{v_{ij}}$ are coordinate unit vectors. The target
is narrower than the source and is shifted in one horizontal and one
vertical coordinate. This choice lets us see whether SRCFM learns the
contraction as well as the shift. In particular, the vertical
displacement must come from the horizontal motion through
$\dot v_{12}=-h_1\dot h_2$.

We draw 16 independent pairs $(w,z)\sim\mu\otimes\nu$ at each update.
For each pair, we form $q=w^{-1}*z$ using 
\eqref{e-groupop} and \eqref{e-groupinv}
and solve
\eqref{e-Xequation} for the skew-symmetric matrix $A$ and
initial control $v_0$. We use nonlinear least squares with up to eight
initial guesses and 1,000 function evaluations per guess. We check
each endpoint using the path formula and a separate integration, with
error limits $10^{-7}(1+\max\{\|w\|_2,\|z\|_2\})$ in the original
coordinates and $10^{-7}(1+\|q\|_2)$ in the relative coordinates.
The independent integration uses relative and absolute tolerances
$10^{-10}$ and $10^{-12}$. If the solver fails these checks, training
stops and saves the pair for replay. We sample 16 independent uniform
times on each computed path and fit its four control components by
mean squared error.

The network takes time and all ten coordinates as inputs. It applies
$\operatorname{arsinh}$ to each state coordinate and uses three
residual blocks of width 128, SiLU activations, and 32 sinusoidal time
features. Denote the raw output by $b_\theta\in\mathbb R^m$. We use the
following expression for the control
\begin{equation}
\label{eq:num-growth-control}
 a_\theta(t,x)=S(x)
       \frac{b_\theta(t,x)}{\sqrt{1+\|b_\theta(t,x)\|^2/C^2}},
 \qquad
 S(h,v)=\left(1+\|h\|^2+\sqrt{1+\|v\|^2}\right)^{1/2}.
\end{equation}
We set $C=20$, which bounds the control norm by $CS(x)$. We train on a
Tesla V100 GPU with 32 GB of memory using AdamW, weight decay
$10^{-6}$, and gradient-norm clipping at 1. The learning rate increases
to $2\times10^{-4}$ over 500 updates and follows a cosine schedule
with $2\times10^{-5}$ at update 200,000. We use the saved parameters
at update 13,406, averaged during training with exponential
decay 0.999.

We integrate the learned field from 512 fresh source points by adaptive
fourth-order Runge--Kutta with step doubling. We save the states on
200 intervals with relative and absolute tolerances $10^{-5}$ and
$10^{-7}$, and repeat on 400 intervals with tolerances
$6.25\times10^{-7}$ and $6.25\times10^{-9}$. All 512 trajectories reach
$t=1$. We use the finer run for the plots and the distribution errors.

\begin{figure}[!htbp]
\centering
\includegraphics[width=0.65\textwidth]{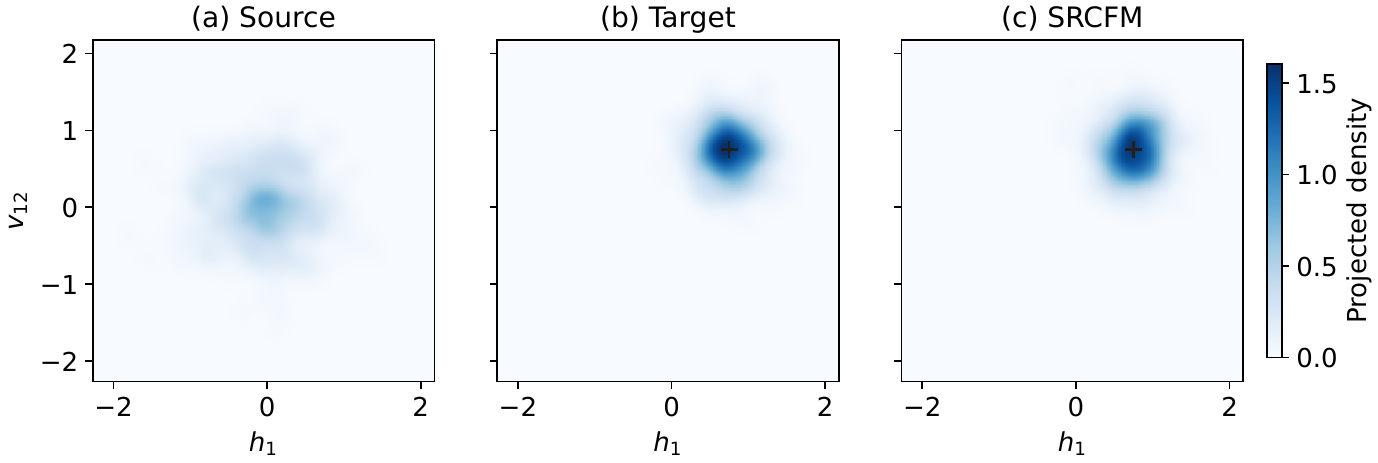}\\[3pt]
\includegraphics[width=0.65\textwidth]{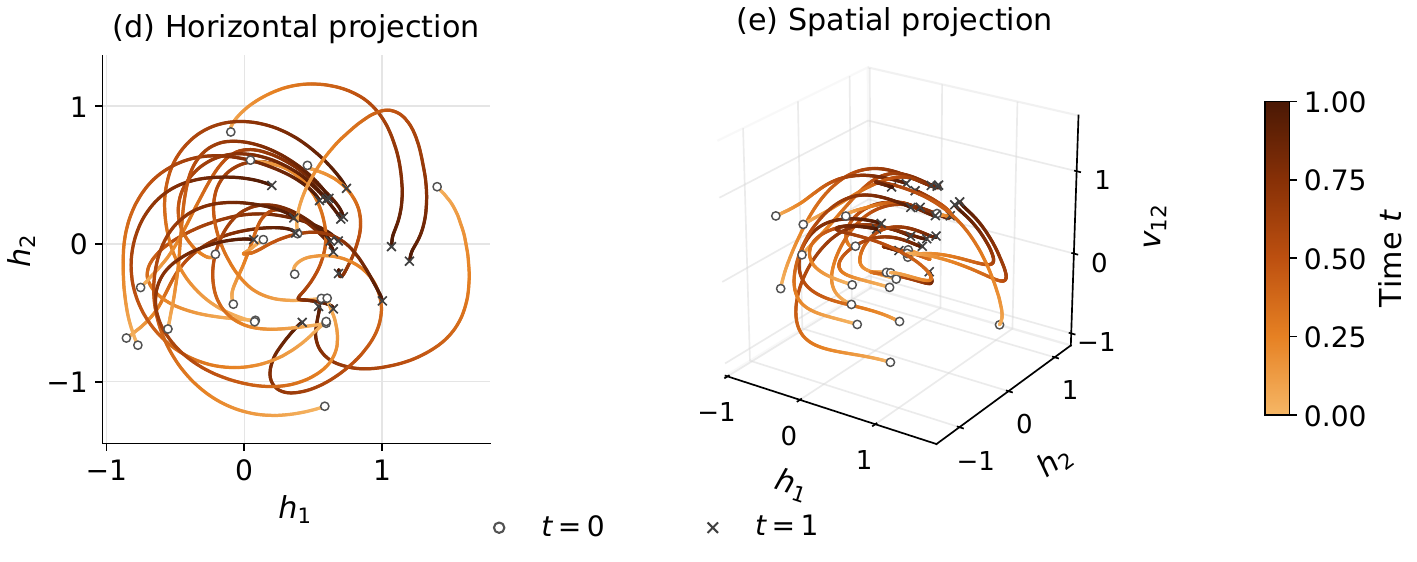}
\caption{Gaussian transport on $\mathbb{R}^{10}$ with the free group structure. Panels (a)--(c) show the
projected source, target, and generated densities using 512 points
each and a common color scale. Panels (d) and (e) show the first 20
generated trajectories at all 401 saved times. Circles and crosses
mark the states at $t=0$ and $t=1$.}
\label{fig:num-gaussian-general}
\end{figure}

Panels (a)--(c) of Figure~\ref{fig:num-gaussian-general} show kernel density
estimates in the $(h_1,v_{12})$ plane. Darker blue indicates a higher
projected density, and the plus signs mark the target mean. The source
density is spread around the origin, while the generated density is
concentrated near the target mean with a spread close to that of the
target. Panels (d) and (e) show the same 20 learned trajectories,
first in the horizontal plane $(h_1,h_2)$ and then with $v_{12}$
included. Following the color from light to dark shows that many of the
paths turn in the horizontal plane before arriving near the shifted
region. Panel (e) shows the change in height that accompanies this
motion, as determined by the horizontal lift.

To check agreement across all ten coordinates, we measure sliced
Wasserstein distance to a reference set of 4,096 target points.
We use 256 fixed random projections and include a separate
512-point target sample as a finite-sample comparison. The distance
drops from $0.40$ for the source to $0.029$ for the generated samples,
compared with $0.023$ for the target sample. Using the same finite
differences, the mean absolute error in the vertical equations
\eqref{e-lift} is about $10^{-5}$. The endpoints change little under
refinement.

\subsection{\texorpdfstring{Endpoint matching on subspaces of $\mathbb{F}_{100,2}=\mathbb{R}^{5050}$}{Endpoint matching on subspaces of F100,2=R5050}}
\label{sec:num-m100-endpoint}

We next check the endpoint accuracy of the piecewise formula in
Theorem~\ref{thm:single-excursion-connector} when the number of controls
is $m=100$. We retain the vertical coordinates $v_{ij}$ indexed by a
set $\mathcal I$ of pairs $1\leq i<j\leq m$, giving state dimension
$d=m+|\mathcal I|$. We take $|\mathcal I|=300$ and $684$ to obtain
$d=400$ and $784$, respectively. Each set includes the 99 consecutive
pairs $(i,i+1)$. We draw the remaining pairs uniformly without
replacement. For each geometry we sample
100 independent pairs of endpoints from
$\mathcal N(0,0.35^2 I_d)\otimes\mathcal N(0,0.35^2 I_d)$.

Given endpoints $w,z$, we again compute $q=w^{-1}*z$ using 
\eqref{e-groupop} and \eqref{e-groupinv}. For the pairs in $\mathcal I$, we use
the entries of the area matrix in \eqref{eq:single-excursion-area}.
We set the other upper-triangular entries to zero, corresponding to the choice
$q^{(i,j)}=-q^i q^j/2$ for the omitted coordinates of the relative
endpoint. After evaluating the piecewise formula, we keep all
horizontal coordinates and the vertical coordinates indexed by
$\mathcal I$. The resulting path satisfies the horizontal lift
\eqref{e-lift} on these retained pairs.

We compare the endpoint from the formula with an independent
computation of the lift. We integrate the vertical velocities along
the loop using 32-node Gauss--Legendre quadrature, first dividing the
loop interval into 16 panels and then into 32 panels. In each case we
measure $\|\gamma(1)-z\|_2$ in all $d$ coordinates. The tolerance for
the two integrated endpoint errors and their endpoint difference is
$10^{-9}(1+\max\{\|w\|_2,\|z\|_2\})$.

\begin{table}[htbp]
\centering\small
\caption{Maximum endpoint errors over 100 pairs for each state
dimension with $m=100$. The last column uses 32 quadrature panels on
the loop.}
\label{tab:num-m100-endpoint}
\begin{tabular}{rrrr}
\toprule
$d$ & Vertical coordinates & Formula error & Integrated error\\
\midrule
400 & 300 & $3.71\times10^{-14}$ & $5.75\times10^{-14}$\\
784 & 684 & $5.97\times10^{-14}$ & $1.17\times10^{-13}$\\
\bottomrule
\end{tabular}
\end{table}

Both quadrature calculations pass the endpoint checks for all pairs.
The errors in Table~\ref{tab:num-m100-endpoint} show that the loop
provides the prescribed vertical displacement to numerical precision.
We use the same construction for training paths in the image example.

\subsection{Comparison of paths between fixed endpoints}

\begin{figure}[h!]
\centering
\includegraphics[width=0.8\textwidth]{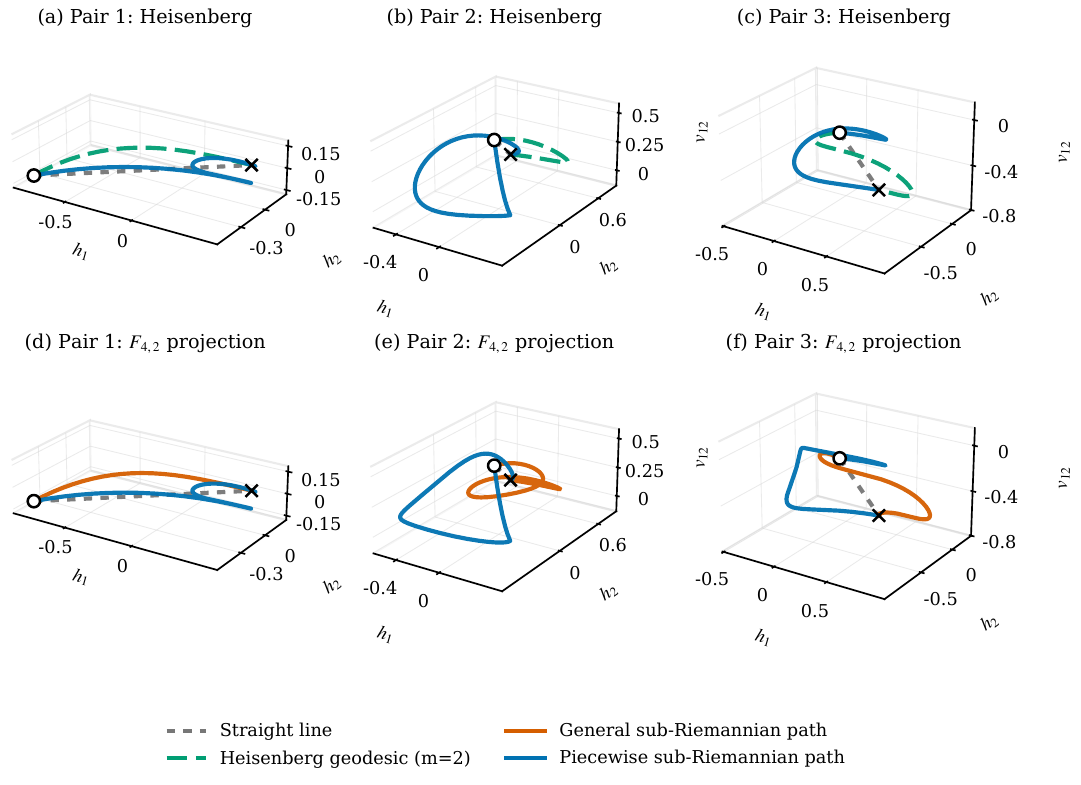}
\caption{Paths between prescribed endpoints. The top row shows
the Heisenberg group, with the scalar geodesic and the $m=2$ general
formula coinciding in green. The bottom row shows paths in $\mathbb{R}^{10}$
projected onto $(h_1,h_2,v_{12})$. Orange denotes the general
sub-Riemannian path, blue the piecewise sub-Riemannian path, and
dashed gray the straight segment. Circles mark sources and crosses
mark targets.}
\label{fig:num-single-excursion-paths}
\end{figure}

We next compare the paths used to construct training targets.
Figure~\ref{fig:num-single-excursion-paths} shows three chosen pairs of
endpoints, one pair in each column. The top row uses the Heisenberg
group. For these pairs, the scalar geodesic calculation agrees with
the general formula for $m=2$, so the two paths appear as one green
curve. The bottom row shows paths computed in all ten coordinates of
$\mathbb{R}^{10}$ and projected onto $(h_1,h_2,v_{12})$.

In each panel the curves start at the same circle and end at the same
cross. In the first column, several of the curves take a similar route in the
projection. In the second column, the piecewise path (blue) makes a
broader turn than the other paths, and the third column shows the
horizontal paths bending to different sides of the straight segment. The piecewise
construction first moves to the target horizontal position and then
uses a loop in the horizontal coordinates to supply the remaining
vertical displacement. The dashed straight segment provides a visual
reference for the difference in shape. For the general paths, we use
at most 12 initial guesses, 300 function evaluations per guess, and
absolute and relative endpoint tolerances of $10^{-9}$. The displayed
general paths reach their endpoints to within $7\times10^{-11}$ in
the full state.

\subsection{SRCFM on \texorpdfstring{$14\times14$}{14x14} MNIST}
\label{sec:num-mnist-single-excursion}

For the final example we use SRCFM to evolve Gaussian samples in the space of resized MNIST images. We resize the images to $14\times14$ by antialiased bilinear interpolation, then rescale intensities to $[-1,1]$. We use 55,000 images for fitting and 5,000 for validation, all from the MNIST training set. The source and target laws are
\[
 \mu=\mathcal N(0,0.5^2 I_{196}),\qquad
 \nu=\frac1{55000}\sum_{a=1}^{55000}\mathcal N(y_a,0.01^2 I_{196}),
\]
where $y_a$ denotes the fitting images in the coordinate order specified below. The target therefore samples a training image uniformly and adds independent Gaussian noise of standard deviation $0.01$ to each pixel. Both laws have densities on $\mathbb R^{196}$.

We use the geometry of Section~\ref{sec:num-m100-endpoint} with $m=100$ and $|\mathcal I|=96$. The 100 horizontal pixels are at the intersections of rows and columns with zero-based indices $\{0,1,3,4,6,7,9,10,12,13\}$. We arrange these coordinates in four chains of 25. Each neighboring pair in a chain contributes one pair to $\mathcal I$. The other 96 pixels are the vertical coordinates, ordered by row. At each update we sample 16 independent pairs from $\mu\otimes\nu$ and evaluate the piecewise path at eight uniform times per pair, giving 128 training examples. The horizontal projection of each path follows a straight segment and then a loop. Its lift supplies the prescribed vertical displacement, and the path has constant horizontal speed almost everywhere.

The network takes time and all 196 coordinates as input and returns the 100 controls $a_\theta$. The sampling equation is
\[
 \dot h_j=a_{\theta,j}(t,h,v),\qquad
 \dot v_{ij}=-h_i a_{\theta,j}(t,h,v),\quad (i,j)\in\mathcal I.
\]
Thus the vertical pixel velocities follow from the controls and the current horizontal coordinates. We use a width-512 network with four residual blocks, SiLU activations, and a 64-dimensional sinusoidal time embedding, for a total of 2,566,372 parameters. The input transformation is $\operatorname{arsinh}((x-c)/s)$, applied coordinatewise with scales at least $0.05$. We use 256 fitting pairs and eight times per pair to fix the centers $c$, the scales $s$, the output scale, the loss scale, and the constant $C=75.4945$ in \eqref{eq:num-growth-control}. The training objective is the mean squared error in the controls divided by the square of the fixed loss scale.

We train on a GPU using AdamW with weight decay $10^{-6}$ and gradient norm capped at 1. The learning rate increases to $2\times10^{-4}$ over the first 1,000 updates, follows a cosine schedule to $2\times10^{-5}$ at update $10^6$, and remains there through update $2\times10^6$. The exponential moving average of the parameters uses decay $0.999$.

\begin{figure}[!htbp]
\centering
\includegraphics[width=0.76\textwidth]{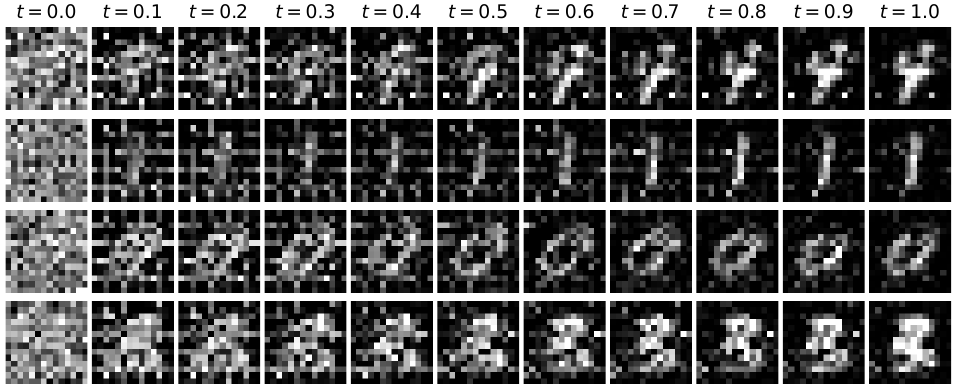}
\caption{Four examples of MNIST generation at $t=0,0.1,\ldots,1$.
Each row follows a different Gaussian input. These examples were
selected from 16 trajectories for their recognizable output shapes.
Black and white correspond to $-1$ and $1$, with clipping applied only
for display.}
\label{fig:num-mnist-single-evolution}
\end{figure}

We integrate the final averaged network from 16 fixed Gaussian inputs
using adaptive RK4 with relative and absolute tolerances of $10^{-4}$ and
$10^{-6}$. Figure~\ref{fig:num-mnist-single-evolution} shows the evolution
from Gaussian inputs on the left to generated images on the right. We
display each state in the original pixel order. As the background
darkens, localized strokes become visible, including a narrow upright
stroke and a loop. Some shapes remain incomplete, and some intensities
fall outside the displayed range.

SRCFM generates each 196-pixel image by directly controlling the velocities of 100 horizontal pixels and obtaining the other 96 velocities from the horizontal lift.  These selected trajectories give a qualitative view of the generation process. Quantitative comparison with the target distribution would require independent generated samples.

\begin{figure}[!htbp]
\centering
\includegraphics[width=0.75\textwidth]{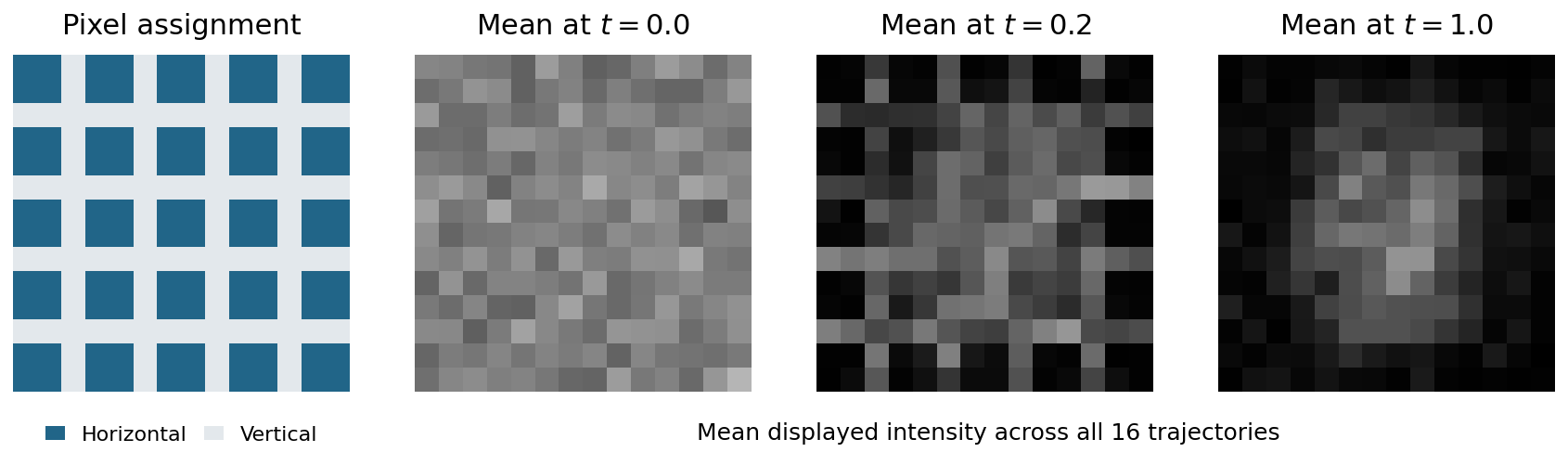}
\caption{Pixel assignment and mean images during MNIST generation.
The left panel shows horizontal pixels in blue and vertical pixels
in light gray. The other panels average the displayed pixel
intensities across all 16 trajectories at $t=0$, $0.2$, and $1$.
Each image is clipped to $[-1,1]$ before averaging.}
\label{fig:num-mnist-coordinate-grid}
\end{figure}

The early columns of Figure~\ref{fig:num-mnist-single-evolution} show a grid. In Figure~\ref{fig:num-mnist-coordinate-grid}, the pixel assignment shows horizontal $2\times2$ blocks separated by rows and columns of vertical pixels. The mean image has no clear grid at $t=0$.  At $t=0.2$, the blocks are darker than the intervening lines. This contrast is present before clipping and weakens by $t=1$.

The two types of pixels obey different velocity equations. The network supplies the horizontal velocity $\dot h_j=a_{\theta,j}$ directly, and the vertical velocity $\dot v_{ij}=-h_i a_{\theta,j}$ also depends on the current horizontal coordinate $h_i$. In this experiment the horizontal pixels darken sooner. Because the two types of pixels are arranged regularly, the difference in their evolution appears as a grid.

\section{Conclusion}
\label{sec:conclusion}

We developed a common training procedure for sub-Riemannian conditional
flow matching using three constructions of horizontal paths. These are
geodesics in the Heisenberg group and its products, paths obtained from
the endpoint equations for free step-two Carnot groups, and explicit
piecewise horizontal paths. The procedure learns the horizontal controls,
which determine the full velocity and enforce the prescribed motion
constraints during generation.

The explicit piecewise construction gives the path and its control
directly from the sampled endpoints, making repeated path evaluation
practical in our higher-dimensional examples. This provides a basis
for extending the method to problems with more controls and state
coordinates. Future work can build on these results by developing
faster path constructions and more robust solvers for the general
endpoint equations.

 \appendix

 \section{More horizontal calculus}
 
    \subsection{Gradients}

    Recall that, for a differentiable function $f:\mathbb{R}^d \to \mathbb{R}$, the directional derivative of $f$ at $x \in \mathbb{R}^d$ in the direction of a vector $w \in \mathbb{R}^d$ is defined as
    $
    \nabla f \cdot w
    $.
    With this in mind, we can write
    \[
    \nabla f 
    = \sum_{j=1}^d \left(\nabla f\cdot \frac{\partial}{\partial x^j}\right) \frac{\partial}{\partial x^j}
    \]
    In other words, we obtain $\nabla f$ by projecting it along each Euclidean coordinate vector.
    Projecting $\nabla f$ instead along horizontal coordinate vectors gives us the following.
    \begin{definition}[Horizontal gradient]    
    The {\em horizontal gradient} of a differentiable function $f:\mathbb{R}^d \to \mathbb{R}$ is defined as
    \[
    \nabla_Hf 
    =
    \sum_{i=1}^m \left(\nabla f \cdot v_i\right) v_i.
    \]
    \end{definition}

    \subsection{Divergence}
    The standard Euclidean divergence of a vector field $u = \sum_{j=1}^d u^j \frac{\partial}{\partial x^j}$ is
    \[
    \operatorname{div}(u) 
    = \sum_{j=1}^d \left( \nabla u^j \cdot \frac{\partial}{\partial x^j} \right).
    \]
    That is, we sum the $j$th Euclidean coordinates of $\nabla u^j$.
    Analogously summing the $i$th horizontal coordinates of $\nabla_H u_i$ gives the following.
    
    \begin{definition}[Horizontal divergence]
    The {\em horizontal divergence} of a horizontal vector field $u = \sum_{i=1}^m u_i v_i$ is defined as
    \[
    \operatorname{div}_H(u) 
    =
    \sum_{i=1}^m \left( \nabla u_i \cdot v_i \right).
    \]
    \end{definition}
    Note that, for a horizontal vector field $u = \sum_{i=1}^m u_i v_i$, we have
    \begin{align*}
        \operatorname{div}(u)
        =
        \operatorname{div}\left( \sum_{i=1}^m u_i v_i \right)
        =
        \sum_{i=1}^m\operatorname{div}\left(  u_i v_i \right)
        =
        \sum_{i=1}^m (\nabla u_i \cdot v_i + u_i \operatorname{div} (v_i))
        =
        \operatorname{div}_H(u) + \sum_{i=1}^m u_i \operatorname{div} (v_i).
    \end{align*}
    This leads to the following well-known fact:
    \begin{proposition}
    If $\operatorname{div}(v_i) = 0$ for $i=1,\dots,m$, then $\operatorname{div}(u) = \operatorname{div}_H(u)$ for any horizontal vector field $u$.
    \end{proposition}

\bibliographystyle{unsrtnat}
\bibliography{ref}

\end{document}